\pdfoutput=1
\RequirePackage{ifpdf}
\ifpdf 
\documentclass[pdftex]{sigma}
\else
\documentclass{sigma}
\fi

\numberwithin{equation}{section}

\usepackage{amssymb,stmaryrd,setspace,bbm}

\usepackage{array}
\newcolumntype{C}{>{$}c<{$}} 

\usepackage{enumitem}
\setitemize{leftmargin=*} 
\setenumerate{leftmargin=*, 
 label=\textnormal{\arabic*.\@}
	}

\usepackage{mathtools} 

\usepackage{tikz}
\usetikzlibrary{calc}
\tikzset{%
	>=latex,
	wt/.style={circle, draw=black, fill=black, inner sep=2pt, outer sep=0pt, minimum size=5pt}, 
	nom/.style={circle, draw=black!20, fill=black!20, inner sep=1pt},
 owt/.style={circle, draw=black, inner sep=2pt, outer sep=0pt, minimum size=5pt},
}

\newcommand{\pd}{\partial} 
\newcommand{\wun}{\mathbbm{1}} 

\newcommand{\bchi}{\bar{\chi}{}} 
\newcommand{\bpsi}{\bar{\psi}{}} 
\newcommand{\bG}{\bar{G}{}}
\newcommand{\bm}{\bar{m}}

\newcommand{\cc}{\mathsf{c}} 
\newcommand{\ee}{\mathrm{e}} 
\newcommand{\ii}{\mathrm{i}} 
\newcommand{\kk}{\mathsf{k}} 
\newcommand{\qq}{\mathsf{q}} 
\newcommand{\yy}{\mathsf{y}} 
\newcommand{\zz}{\mathsf{z}} 

\newcommand{\NS}{\textnormal{NS}}
\newcommand{\R}{\textnormal{R}}

\renewcommand{\ge}{\geqslant} 
\renewcommand{\le}{\leqslant} 

\DeclarePairedDelimiter{\brac}{\lparen}{\rparen} 
\DeclarePairedDelimiter{\sqbrac}{\lbrack}{\rbrack} 
\DeclarePairedDelimiter{\dbrac}{\llbracket}{\rrbracket} 

\DeclarePairedDelimiterX{\comm}[2]{\lbrack}{\rbrack}{#1 , #2} 
\DeclarePairedDelimiterX{\acomm}[2]{\lbrack}{\rbrack}{#1 , #2} 
\DeclarePairedDelimiterX{\bilin}[2]{\langle}{\rangle}{#1 , #2} 

\newcommand{\no}[1]{\mathopen{:} #1 \mathclose{:}} 
\DeclarePairedDelimiter{\ket}{\lvert}{\rangle}
\DeclarePairedDelimiterXPP{\nsket}[1]{}{\lvert}{\rangle}{^{\NS}}{#1}
\DeclarePairedDelimiterXPP{\rket}[1]{}{\lvert}{\rangle}{^{\R}}{#1}

\DeclareMathOperator{\spn}{span}

\newcommand{\blank}{{-}} 

\newcommand{\Ra}{\Rightarrow}

\newcommand{\lra}{\longrightarrow}
\newcommand{\ira}{\hookrightarrow}

\newcommand{\ses}[3]{0 \lra #1 \lra #2 \lra #3 \lra 0}

\newcommand{\fld}[1]{\mathbb{#1}} 
\newcommand{\alg}[1]{\mathfrak{#1}} 
\newcommand{\grp}[1]{\mathsf{#1}} 
\newcommand{\VOA}[1]{\mathsf{#1}} 
\newcommand{\Mod}[1]{\mathcal{#1}} 
\newcommand{\categ}[1]{\mathrm{#1}} 

\newcommand{\ZZ}{\fld{Z}}

\newcommand{\QQ}{\fld{Q}}

\newcommand{\CC}{\fld{C}}

\newcommand{\affine}[1]{\widehat{#1}}
\newcommand{\SLA}[2]{\alg{#1}_{#2}} 
\newcommand{\SLSA}[3]{\alg{#1}_{#2 \vert #3}} 
\newcommand{\AKMSA}[3]{\affine{\alg{#1}}_{#2 \vert #3}} 

\newcommand{\gltwo}{\SLA{gl}{2}}
\newcommand{\sltwo}{\SLA{sl}{2}}

\newcommand{\slthree}{\SLA{sl}{3}}
\newcommand{\pslone}{\SLSA{psl}{1}{1}}
\newcommand{\psl}{\SLSA{psl}{2}{2}}
\newcommand{\apsl}{\AKMSA{psl}{2}{2}}

\DeclarePairedDelimiterXPP{\envalg}[1]{\mathsf{U}}{\lparen}{\rparen}{}{#1} 

\newcommand{\bgsymb}{\VOA{BG}}
\newcommand{\fgsymb}{\VOA{FG}}
\newcommand{\sfsymb}{\VOA{SF}}
\newcommand{\lsymb}{\Pi}

\newcommand{\bgvoa}{\bgsymb} 
\newcommand{\fgvoa}{\fgsymb} 
\newcommand{\sfvoa}{\sfsymb} 
\newcommand{\lvoa}{\lsymb} 

\newcommand{\uaff}[2]{\VOA{V}^{#1}(#2)} 
\newcommand{\saff}[2]{\VOA{L}_{#1}(#2)} 
\newcommand{\upsl}[1][\kk]{\uaff{#1}{\psl}} 
\newcommand{\usl}[1][\kk]{\uaff{#1}{\sltwo}} 
\newcommand{\ssl}[1][\kk]{\saff{#1}{\sltwo}} 

\newcommand{\uwalg}[2]{\VOA{W}^{#1}_{\mathrm{#2}}} 
\newcommand{\swalg}[2]{\VOA{W}_{#1}^{\mathrm{#2}}} 
\newcommand{\upr}[1][\kk]{\uwalg{#1}{pr}} 
\newcommand{\umin}[1][\kk]{\uwalg{#1}{min}} 
\newcommand{\spr}[1][\kk]{\swalg{#1}{pr}} 
\newcommand{\smin}[1][\kk]{\swalg{#1}{min}} 
\newcommand{\sminpmhalf}{\smin[\pm 1/2]}

\newcommand{\sfnsvac}{\ket{\NS}} 
\newcommand{\sfrvac}{\ket{\R}} 

\newcommand{\nconjsymb}{\mathsf{w}} 
\newcommand{\nconj}{\nconjsymb \brac}
\newcommand{\nsfsymb}[1]{\sigma^{#1}} 
\newcommand{\nsf}[1]{\nsfsymb{#1} \brac}

\DeclarePairedDelimiterXPP{\zhu}[1]{\mathsf{Zhu}}{\big\lbrack}{\big\rbrack}{}{#1} 
\newcommand{\zh}[1]{#1_{\circ}} 

\newcommand{\zhver}[1]{\zh{\Mod{V}}^{#1}} 

\newcommand{\lmod}[1]{\lsymb_{[#1]}} 

\newcommand{\sfnsmod}{\Mod{S}^{\NS}} 
\newcommand{\sfrmod}{\Mod{S}^{\R}} 

\newcommand{\nnsind}[1]{\Mod{N}^{\NS}_{[#1]}} 
\newcommand{\nnsrel}[1]{\Mod{M}^{\NS}_{\vphantom{[#1]}\smash{\dbrac{#1}}}} 
\newcommand{\nnsihw}[1]{\Mod{L}^{\NS}_{#1}} 
\newcommand{\nrind}[1]{\Mod{N}^{\R}_{[#1]}}
\newcommand{\nrrel}[1]{\Mod{M}^{\R}_{\vphantom{[#1]}\smash{\dbrac{#1}}}}
\newcommand{\nrihw}[1]{\Mod{L}^{\R}_{#1}}

\newcommand{\nrverma}[1]{\Mod{W}^{\R}_{[#1]}}
\newcommand{\nrver}[1]{\Mod{V}^{\R}_{\vphantom{[#1]}\smash{\dbrac{#1}}}}
\newcommand{\nrlog}[1]{\Mod{Q}^{\R}_{[#1]}}
\newcommand{\nrsublog}[1]{\Mod{P}^{\R}_{\vphantom{[#1]}\smash{\dbrac{#1}}}}

\newcommand{\adfunc}[1]{\categ{Ad}_{[#1]}}
\DeclareMathOperator{\res}{\categ{Res}}

\DeclareMathOperator{\tr}{tr}
\newcommand{\traceover}[1]{\tr_{\raisebox{-2pt}{$\scriptstyle #1$}}} 

\newcommand{\Gr}[1]{\sqbrac[\big]{#1}} 

\DeclareMathOperator{\chmap}{ch}

\newcommand{\ch}[1]{\chmap^+ \Gr{#1}} 
\newcommand{\fch}[2]{\ch{#1} \brac[\big]{#2}} 
\newcommand{\sch}[1]{\chmap^- \Gr{#1}} 
\newcommand{\fsch}[2]{\sch{#1} \brac[\big]{#2}}
\newcommand{\gch}[1]{\chmap^{\pm} \Gr{#1}} 
\newcommand{\fgch}[2]{\gch{#1} \brac[\big]{#2}}

\newcommand{\fms}{Friedan--Martinec--Shenker}

\newcommand{\km}{Kac--Moody}

\newcommand{\ns}{Neveu--Schwarz}
\newcommand{\pbw}{Poincar\'{e}--Birkhoff--Witt}
\newcommand{\wzw}{Wess--Zumino--Witten}

\newcommand{\fdim}{finite-dimensional}

\newcommand{\lhs}{left-hand side}
\newcommand{\lhss}{\lhs s}
\newcommand{\rhs}{right-hand side}
\newcommand{\rhss}{\rhs s}

\newcommand{\hw}{highest-weight}
\newcommand{\hwv}{\hw\ vector}
\newcommand{\hwvs}{\hwv s}
\newcommand{\hwm}{\hw\ module}
\newcommand{\hwms}{\hwm s}

\newcommand{\rhw}{relaxed highest-weight}

\newcommand{\rhwm}{\rhw\ module}
\newcommand{\rhwms}{\rhwm s}

\newcommand{\chw}{conjugate highest-weight}
\newcommand{\chwv}{\chw\ vector}

\newcommand{\vo}{vertex-operator}
\newcommand{\voa}{\vo\ algebra}
\newcommand{\voas}{\voa s}
\newcommand{\svoa}{\vo\ superalgebra}
\newcommand{\svoas}{\svoa s}

\newcommand{\va}{vertex algebra}
\newcommand{\vas}{\va s}
\newcommand{\sva}{vertex superalgebra}
\newcommand{\svas}{\sva s}
\newcommand{\vsa}{vertex subalgebra}

\newcommand{\emt}{energy-momentum tensor}

\newcommand{\ope}{operator product expansion}
\newcommand{\opes}{\ope s}
\newcommand{\qhr}{quantum Hamiltonian reduction} 

\newtheorem{Theorem}{Theorem}[section]
\newtheorem{Corollary}[Theorem]{Corollary}
\newtheorem{Lemma}[Theorem]{Lemma}
\newtheorem{Proposition}[Theorem]{Proposition}

\begin{document}

\allowdisplaybreaks

\newcommand{\arXivNumber}{2509.04795}

\renewcommand{\thefootnote}{}

\renewcommand{\PaperNumber}{030}

\FirstPageHeading

\ShortArticleName{The Principal W-Algebra of $\mathfrak{psl}_{2|2}$}

\ArticleName{The Principal W-Algebra of $\boldsymbol{\mathfrak{psl}_{2|2}}$\footnote{This paper is a~contribution to the Special Issue on Recent Advances in Vertex Operator Algebras in honor of James Lepowsky. The~full collection is available at \href{https://sigma-journal.com/Lepowsky.html}{https://sigma-journal.com/Lepowsky.html}}}

\Author{Zachary FEHILY~$^{\rm a}$, Christopher RAYMOND~$^{\rm b}$ and David RIDOUT~$^{\rm a}$}

\AuthorNameForHeading{Z.~Fehily, C.~Raymond and D.~Ridout}

\Address{$^{\rm a)}$~School of Mathematics and Statistics, University of Melbourne, Australia}
\EmailD{\mail{zacfehily@gmail.com}, \mail{david.ridout@unimelb.edu.au}}

\Address{$^{\rm b)}$~Department of Mathematics, University of Hamburg, Germany}
\EmailD{\mail{christopher.raymond@uni-hamburg.de}}

\ArticleDates{Received September 08, 2025, in final form February 18, 2026; Published online March 25, 2026}

\Abstract{We study the structure and representation theory of the principal W-algebra~\smash{$\mathsf{W}_{\rm pr}^{\mathsf{k}}$} of~$\mathsf{V}^{\mathsf{k}}(\mathfrak{psl}_{2|2})$. The defining operator product expansions are computed, as is the Zhu algebra, and these results are used to classify irreducible highest-weight modules. In particular, for $\mathsf{k} = \pm \frac{1}{2}$, \smash{$\mathsf{W}_{\rm pr}^{\mathsf{k}}$} is not simple and the corresponding simple quotient is the symplectic fermion vertex algebra. We use this fact, along with inverse Hamiltonian reduction, to study relaxed highest-weight and logarithmic modules for the small $N=4$ superconformal algebra at central charges $-9$ and $-3$.}

\Keywords{vertex-operator algebras; conformal field theory; representation theory}

\Classification{17B67; 17B69; 81T40}

\renewcommand{\thefootnote}{\arabic{footnote}}
\setcounter{footnote}{0}

\section{Introduction}

The simple Lie superalgebra \smash{$\psl$} is somewhat unusual among the basic classical examples in that it has odd roots of multiplicity $2$, which may be regarded as being both positive and negative.
Nevertheless, it has received an enormous amount of attention due to its applications in physics.
In particular, the ``baby'' version of the AdS/CFT correspondence on $\grp{AdS}_3 \times \grp{S}^3$ supersymmetrises to a sigma model on the Lie supergroup \smash{$\grp{PSU}_{1,1\vert2}$} \cite{BerCon99}, leading to interest in \wzw\ models with \smash{$\apsl$} symmetry (see, for example, \cite{EbeWor18,FerL1p24,GabMas11,GabU2223,GotWZN06,TroMas11}).

Intertwined with this story is that of the (small) $N=4$ superconformal algebra.
Introduced initially to describe spacetime supersymmetries of heterotic strings, its unitary representations were first studied in \cite{EguUni87}, motivated by applications to string theory on hyperK\"{a}hler manifolds.
More recently, there has been a resurgence in interest in $N=4$ representation theory because of its role in such things as Mathieu moonshine and its variants \cite{EguNot10}.
However, this representation theory is still poorly understood in general.
Classifications are only known for the unitary \hwms\ \cite{KacUni22,KacSpe25} and the irreducibles in the case of central charge $\cc=-9$ \cite{AdaRea14}.
The $\cc= -9$ algebra is also the $p=2$ member of the family of extensions $\mathcal{V}^{(p)}$ of $\ssl$ at~${\kk=-2+\frac{1}{p}}$, introduced in \cite{AdaRea14} and studied further in \cite{AdaVer20}.

In a different direction, recent developments in the $3$d/$2$d (and $4$d/$2$d) correspondence has renewed interest in theories with $N=4$ supersymmetry.
In \cite{CosVer18}, it was argued that generic $3$-dimensional gauge theories with an appropriate boundary condition support \voas\ on the boundary.
The motivation here was that the data of the \voa, for example, its fusion algebra, could be used to investigate a related topological quantum field theory (analogous to the well-known story connecting \wzw\ models and $3$d Chern--Simons theories).

A rich source of $3$-dimensional gauge theories are the quiver gauge theories \cite{DouDbr96}, in which the field content is determined by specifying a quiver.
These theories have an associated geometry described by the corresponding quiver varieties of Nakajima \cite{NakIns94,NakQui98}.
One example, analysed in~\cite{AraHil23}, is associated with the Jordan quiver for which the associated quiver variety is the Hilbert scheme $\operatorname{Hilb}^{n}(\CC)$ of $n$-points in the complex plane.
The authors use chiral quantisation to construct a \voa\ for each $n$ that has $\operatorname{Hilb}^{n}(\CC)$ as its associated variety.
They also showed that for $n=2$, their \voa\ is the small $N=4$ superconformal algebra at $\cc = -9$ (tensored with some additional free-field vertex algebras).
As such, a detailed understanding of $N=4$ representation theory, its modular data and fusion rules, will shed light on the associated topological quantum field theory.

With this goal in mind, our interest is the relationship between the \svas\ constructed from \smash{$\psl$} and $N=4$.
In the formalism developed by Kac, Roan and Wakimoto, these are related by \qhr: the small $N=4$ superconformal algebra is a minimal reduction of the affine \sva\ constructed from \smash{$\psl$} \cite{KacQua03,KacQua03b}.
But there is also a principal reduction which appears to have received little to no attention in the literature.
This paper then amounts to the first steps towards understanding this reduction and its representation theory.

Given the recent interest in $N=4$ representation theory, we aim to connect this to the representation theory of the principal reduction.
Happily, the formalism of inverse quantum Hamiltonian reduction is designed for this purpose.
Originally introduced by Semikhatov \cite{SemInv94} for the universal Virasoro and $\sltwo$ vertex algebras, inverse reduction was extended to the simple quotients by Adamovi\'{c} in \cite{AdaRea17}.
The latter also studied certain restriction functors that connect the representation theories, specifically sending irreducible \hw\ Virasoro modules to (spectral flows of) \rhw\ $\ssl$-modules.
Subsequent work generalised this to higher ranks and developed methods to prove that these $\ssl$-modules were generically irreducible \cite{AdaRea20} and that every irreducible relaxed $\ssl$-module could be obtained in this fashion~\cite{AdaWei23}.
Other recent works addressing inverse reduction include \cite{AdaNap24,AdaRel21,CreSL225,FasCon24,FasVir24,FasMod24,FehSub21,FehInv23,FehMod21}.

Inverse reduction therefore suggests a general path to understanding the representation theory of a nonrational simple affine \voa\ or W-algebra (associated to a simple Lie superalgebra $\alg{g}$ at level $\kk$).
One starts with the representation theory of the ``exceptional'' W-algebra of level $\kk$ \cite{AraRat19} and then uses inverse reduction \`{a} la Adamovi\'{c} to reconstruct the representation theory of each successive ``less reduced'' W-algebra.
In this way, one ascends through the poset of level-$\kk$ W-algebras corresponding to $\alg{g}$ and $\kk$ (parametrised by the even nilpotent orbits of $\alg{g}$), iteratively building up their representation theories.

We remark that when $\alg{g}$ is a Lie algebra and $\kk$ is (co)admissible, the exceptional W-algebra is rational \cite{AraRat19} and so its representation theory is well understood.
This need not be the case if $\kk$ is nonadmissible, see \cite[Section~5.2]{AdaWei23}, or $\alg{g}$ is not a Lie algebra (or $\SLSA{osp}{1}{2n}$), see \cite{CreSL225}.
In the case of interest here ($\alg{g} = \psl$), we identify the ``exceptional'' W-algebra for $\kk = \pm \frac{1}{2}$ as the symplectic fermions \svoa\ \cite{KauSym00}.
Inverse reduction then allows us to investigate the representation theory of the small $N=4$ superconformal algebra at central charge $\cc=-9$ (and also at $\cc=-3$).

\subsection{Outline and results}

We begin in Section~\ref{sec:principal} by recalling some basic facts about the Lie superalgebra \smash{$\psl$}.
We then compute the principal W-algebra \smash{$\upr$} from the universal affine \vo\ superalgebra~\smash{$\upsl$} via quantum Hamiltonian reduction.
The \opes\ that define~$\upr$ are given in Theorem~\ref{thm:propes} and we record that \smash{$\kk=\pm\frac{1}{2}$} is a collapsing level with simple quotient isomorphic to the symplectic fermions \svoa.
We also determine the levels for which \smash{$\upr$} is not simple (see Theorem~\ref{thm:prsimple}).

In Section~\ref{sec:reps}, we study the lower-bounded weight modules of \smash{$\upr$} by computing the Zhu algebra \smash{$\zhu{\upr}$} and determining the irreducible weight \smash{$\zhu{\upr}$}-modules.
Every irreducible lower-bounded \smash{$\upr$}-module is shown to be highest-weight and the dimensions of their top spaces are determined.
In particular, the top spaces are always \fdim.
This confirms that the representation theory of \smash{$\upr$} is not too dissimilar to those of the principal W-algebras with~$\alg{g}$ a simple Lie algebra.
Given that the symplectic fermions \svoa\ is ``log-rational'', meaning $C_2$-cofinite but not rational, it is reasonable to ask if the simple quotient of~$\upr$ is also log-rational (when the latter is not simple).
We show that this question is interesting by analysing two further examples.

Section~\ref{sec:InvRed} then applies inverse Hamiltonian reduction to \smash{$\upr$}.
We first realise the minimal W-algebra \smash{$\umin$} of \smash{$\psl$} inside \smash{$\upr \otimes \lvoa$}, where \smash{$\lvoa$} is the usual bosonisation of symplectic bosons.
The argument determining when this realisation descends to the simple quotients $\smin$ and $\spr$ is then sketched.
Granting this, Section~\ref{sec:k=+-1/2} restricts to the case $\kk=\pm\frac{1}{2}$ for which inverse reduction translates the well-known representation theory of the symplectic fermions \svoa\ into that of \smash{$\sminpmhalf$}.
We thereby recover the classification of irreducible \hwms\ \cite{AdaRea14} and extend it to irreducible \rhwms\ (in both the \ns\ and Ramond sectors).
The (super)characters of the latter are easily computed, though we leave the potentially subtle analysis of their modularity to future work.

In anticipation of this modularity study, we turn to the degenerations and spectral flows of the relaxed families in Section~\ref{sec:degsfk=+-1/2}.
Interestingly, each member of these families decomposes as a direct sum of two generically irreducible \rhwms.
We also exhibit the structures of the reducible examples using Loewy diagrams.
Here, we observe a qualitative difference between the $\kk=-\frac{1}{2}$ and $\kk=+\frac{1}{2}$ cases: unlike the former, one of the reducible relaxed modules of the latter has composition factors whose minimal conformal weights are strictly greater than that of the top space.
We conclude in Section~\ref{sec:k=+-1/2log} by applying inverse reduction to the logarithmic module of symplectic fermions, thereby obtaining an uncountably infinite number of logarithmic \smash{$\sminpmhalf$}-modules.
The structures of these modules are completely elucidated, even at degeneration points.

\section{Computing the principal W-algebra} \label{sec:principal}

In this section, we introduce the principal W-algebras associated to the simple complex Lie superalgebra \smash{$\psl$}.
For this, we follow the seminal work of Kac, Roan and Wakimoto \cite{KacQua03} to identify generators and compute their defining \opes\ explicitly.

\subsection[The Lie superalgebra psl]{The Lie superalgebra $\boldsymbol{\psl}$} \label{sec:psl}

Recall that \smash{$\SLSA{sl}{2}{2}$} is the Lie superalgebra whose elements are the endomorphisms of the $\ZZ_2$-graded vector superspace \smash{$\CC^{2\vert2}$} with zero supertrace.
The quotient of this superalgebra by the ideal spanned by the identity matrix is denoted by \smash{$\psl$}.
Let $E_{ij}$, $i,j=1,\dots,4$, denote the $4\times4$ matrix whose $(i,j)$-th entry is $1$ and whose other entries are $0$.
We fix a basis of \smash{$\psl$} as follows, indicating parity and understanding that each $E_{ij}$ is implicitly replaced by its image in the quotient\looseness=-1%
\begin{equation}\label{eq:pslbasis}
	\begin{aligned}
	&	\begin{split}
			& E^1 = E_{12}, \quad H^1 = E_{11}-E_{22}, \quad F^1 = E_{21}, \\
			& E^2= E_{34}, \quad H^2 = E_{33}-E_{44}, \quad F^2 = E_{43}
		\end{split}
		& {\rm (even)}, \\
		& \begin{split}
		&	e^{++} = +E_{14}, \quad e^{+-} = -E_{13}, \quad e^{-+} = +E_{24}, \quad e^{--} = -E_{23}, \\
		&	f^{++} = -E_{32}, \quad f^{+-} = -E_{42}, \quad f^{-+} = +E_{31}, \quad f^{--} = +E_{41}
		\end{split}
		& {\rm (odd)}.
	\end{aligned}
\end{equation}

\begin{figure}
	\centering
	\begin{tikzpicture}[<->,scale=1]
		\draw[gray] (-2.5,0) -- (2.5,0);
		\draw[gray] (0,-2.5) -- (0,2.5);
		\node[wt,label=above:{$E^1$}] at (2,0) {};
		\node[wt,label=below:{$F^1$}] at (-2,0) {};
		\node[wt,label=left:{$E^2$}] at (0,2) {};
		\node[wt,label=right:{$F^2$}] at (0,-2) {};
		\node[wt,label=above:{$e^{++}$},label=right:{$f^{++}$}] at (1,1) {};
		\node[wt,label=right:{$e^{+-}$},label=below:{$f^{+-}$}] at (1,-1) {};
		\node[wt,label=left:{$e^{-+}$},label=above:{$f^{-+}$}] at (-1,1) {};
		\node[wt,label=below:{$e^{--}$},label=left:{$f^{--}$}] at (-1,-1) {};
	\end{tikzpicture}
\caption{%
	The roots of the simple Lie superalgebra \smash{$\psl$}, labelled by their root vectors (the odd roots have multiplicity $2$).
	The horizontal axis indicates the eigenvalue under the adjoint action of $H^1$ while the vertical axis records it for $H^2$.
} \label{fig:roots}
\end{figure}
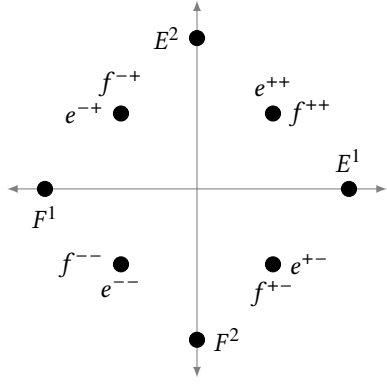

Here are some useful facts about \smash{$\psl$} and the basis that we have chosen:
\begin{itemize}\itemsep=0pt
	\item The even subalgebra is isomorphic to $\sltwo \oplus \sltwo$.
	Both $\big\{ E^1,H^1,F^1\big\}$ and $\big\{ E^2,H^2,F^2\big\}$ are $\sltwo$-bases.
	\item The four odd elements labelled with an $e$ span a copy of the tensor product of two fundamental $\sltwo$-modules.
	The signs in \eqref{eq:pslbasis} ensure that the actions of the $E^i$ and $F^i$, $i=1,2$, just exchange the $i$-th $\pm$ index.
	For example,
$
		F^1 e^{++} = e^{-+}$ and $ E^2 e^{+-} = e^{++}$.
	\item The previous property also holds for the odd elements labelled with an $f$.
	\item The odd roots have multiplicity $2$.
	The corresponding root vectors are nilpotent: $\big\lbrack e^{\pm\pm},e^{\pm\pm}\big\rbrack = \big\lbrack f^{\pm\pm},f^{\pm\pm}\big\rbrack = 0$.
	(Here and throughout, $\comm{\cdot}{\cdot}$ denotes a Lie superbracket, modelled on the anticommutator if both elements are odd and on the commutator if at least one element is even.)
	\item In fact, the odd elements labelled by an $e$ all anticommute with one another.
	The same is true for the odd elements labelled by an $f$.
	\item The nonzero anticommutators are as follows:
		\begin{gather*}
	\acomm[\big]{e^{++}}{f^{+-}} = -E^1, \qquad \acomm[\big]{e^{--}}{f^{++}} = -\tfrac{1}{2}\big(H^1-H^2\big), \qquad \acomm[\big]{e^{--}}{f^{+-}} = +F^2, \\
			\acomm[\big]{e^{+-}}{f^{++}} = +E^1, \qquad \acomm[\big]{e^{+-}}{f^{-+}} = -\tfrac{1}{2}\big(H^1+H^2\big), \qquad \acomm[\big]{e^{+-}}{f^{--}} = -F^2, \\
			\acomm[\big]{e^{-+}}{f^{++}} = -E^2, \qquad \acomm[\big]{e^{-+}}{f^{+-}} = +\tfrac{1}{2}\big(H^1+H^2\big), \qquad \acomm[\big]{e^{-+}}{f^{--}} = +F^1, \\
			\acomm[\big]{e^{++}}{f^{-+}} = +E^2, \qquad\acomm[\big]{e^{++}}{f^{--}} = +\tfrac{1}{2}\big(H^1-H^2\big), \qquad \acomm[\big]{e^{--}}{f^{-+}} = -F^1.
	\end{gather*}
	\item The Killing form is identically zero, but the supertrace form in the defining representation~\eqref{eq:pslbasis} is supersymmetric and nondegenerate.
	The nonzero entries of this supertrace form are as follows:
	\begin{alignat}{5}
			&\big\langle E^1,F^1\big\rangle = +1, \qquad&& \big\langle H^1,H^1\big\rangle = +2, \qquad&& \big\langle E^2,F^2\big\rangle = -1, \qquad&& \big\langle H^2,H^2\big\rangle = -2,&\nonumber \\
			&\big\langle e^{++},f^{--}\big\rangle = +1, \qquad&&\big\langle e^{+-},f^{-+}\big\rangle = -1, \qquad&& \big\langle e^{-+},f^{+-}\big\rangle = -1, \qquad && & \nonumber\\
& \big\langle e^{--},f^{++}\big\rangle = +1. &&&&&& &\label{eq:supertrace}
		\end{alignat}
\end{itemize}
We picture the root system of \smash{$\psl$} in Figure~\ref{fig:roots}, indicating by position the eigenvalues with respect to the basis $\bigl\{H^1,H^2\bigr\}$ of the Cartan subalgebra.

The Weyl group of \smash{$\psl$} is, by definition, that of the even subalgebra $\sltwo \oplus \sltwo$, hence it is isomorphic to $\ZZ_2 \times \ZZ_2$.
However, the automorphism group of the root system is clearly isomorphic to $\grp{D}_4$.
The additional automorphisms are generated by the reflection that fixes the odd root of weight $(1,1)$.
This automorphism corresponds to swapping the two copies of $\sltwo$.
It moreover lifts to an order-$2$ automorphism $\omega$ of \smash{$\psl$} defined by
\begin{alignat*}{5}
		& E^1 \longleftrightarrow E^2, \qquad&& H^1 \longleftrightarrow H^2, \qquad&& F^1 \longleftrightarrow F^2,\qquad && & \\
		&e^{++} \longleftrightarrow f^{++}, \qquad&& e^{+-} \longleftrightarrow f^{-+}, \qquad&& e^{-+} \longleftrightarrow f^{+-}, \qquad&& e^{--} \longleftrightarrow f^{--}.&
	\end{alignat*}
Interestingly, this automorphism does not leave the supertrace form invariant, but instead negates it
\begin{equation} \label{eq:swapsupertrace}
	\bilin[\big]{\omega(A)}{\omega(B)} = -\bilin{A}{B}, \qquad A,B \in \psl.
\end{equation}

\subsection{The principal reduction} \label{sec:Walg}

The (untwisted) affine \km\ superalgebra \smash{$\apsl$} is spanned by elements $A_n$ and $K$, for~${A \in \psl}$ and $n \in \ZZ$, where $K$ is even and $A_n$ inherits its parity from that of $A$.
As usual, $K$ is central and the Lie superbrackets for the $A_n$ are given by
\[
	\comm{A_m}{B_n} = \comm{A}{B}_{m+n} + m \bilin{A}{B} \delta_{m+n,0} K, \qquad A,B \in \psl,\quad m,n \in \ZZ.
\]
Because of \eqref{eq:swapsupertrace}, the automorphism $\omega$ of \smash{$\psl$} lifts to an automorphism of \smash{$\apsl$} by defining
\[
	\affine{\omega}(A_n) = \omega(A)_n \qquad \text{and} \qquad \affine{\omega}(K) = -K, \qquad A \in \psl,\quad n \in \ZZ.
\]
The universal level-$\kk$ affine \sva\ \smash{$\upsl$} is defined, as always, to be strongly and freely generated by the fields \smash{$A(z) = \sum_{n\in\ZZ} A_n z^{-n-1}$}, $A \in \psl$, subject to the \opes
\[
	A(z) B(w) \sim \frac{\bilin{A}{B} \kk \wun}{(z-w)^2} + \frac{\comm{A}{B}(z)}{z-w}, \qquad A,B \in \psl.
\]
Here, $\wun$ denotes the identity field and the $A_n$ furnish a representation of \smash{$\apsl$} on \smash{$\upsl$} in~which $K$ acts as multiplication by the level $\kk \in \CC$.
We note that the automorphism $\affine{\omega}$ of \smash{$\apsl$} again lifts, but this time to an \emph{isomorphism} between \smash{$\upsl$} and $\upsl[-\kk]$.
It follows, in particular, that the representation theory of this \sva\ (and any of its quotients) is independent of the sign of the level.

Because the supertrace form equation~\eqref{eq:supertrace} is nondegenerate, the Sugawara construction equips \smash{$\upsl$} with the structure of a \svoa, when the level is noncritical.
The dual Coxeter number of \smash{$\psl$} is $0$ because the Killing form vanishes, so for all $\kk\ne0$, the \emt\ takes the form
\begin{align*}
	T = {}&\tfrac{1}{2\kk} \bigl( \tfrac{1}{2} \no{H^1H^1} + \no{E^1F^1} + \no{F^1E^1} - \tfrac{1}{2} \no{H^2H^2} - \no{E^2F^2} - \no{F^2E^2} \\
	& - \big(\no{e^{++}f^{--}} - \no{f^{--}e^{++}} - \no{e^{+-}f^{-+}} + \no{f^{-+}e^{+-}} - \no{e^{-+}f^{+-}} \\
	& + \no{f^{+-}e^{-+}} + \no{e^{--}f^{++}} - \no{f^{++}e^{--}}\big) \bigr).
\end{align*}
The central charge is the superdimension of \smash{$\psl$}, namely $-2$.

As there are two nilpotent orbits in $\sltwo$, there are four even orbits in \smash{$\psl$}.
Along with the zero orbit, there are two minimal orbits, for which $F^1$ and $F^2$ are representatives, and the principal one, for which $F^1+F^2$ is a representative.
It is well known that quantum Hamiltonian reduction~$\umin$ with respect to one of the minimal nilpotents is isomorphic to the universal (small) $N=4$ superconformal \svoa\ of central charge $-6(\kk+1)$, see~\cite[Remark~4.1]{KacQua03} and \cite[Section~8.4]{KacQua03b}.
Applying $\affine{\omega}$, we see that the other minimal reduction is isomorphic to $\umin[-\kk]$ with central charge $6(\kk-1)$.
We leave further consideration of these minimal reductions to Section~\ref{sec:InvRed}.

Our interest is in the principal reduction \smash{$\upr$}, which does not seem to have received as much attention.
We observe that $\big\{E^1+E^2,H^1+H^2,F^1+F^2\big\}$ is an $\sltwo$-triple and that the adjoint action of $H^1+H^2$ gives a good even grading on \smash{$\psl$}
\[
		\psl = \psl^{(-2)} \oplus \psl^{(0)} \oplus \psl^{(+2)}, \qquad
	\begin{aligned}
		&\psl^{(+2)} = \spn\bigl\{E^1,E^2 \mid e^{++},f^{++}\bigr\}, \\
		&\psl^{(0)} = \spn\bigl\{H^1,H^2 \mid e^{+-},e^{-+},f^{+-},f^{-+}\bigr\}, \\
		&\psl^{(-2)} = \spn\bigl\{F^1,F^2 \mid e^{--},f^{--}\bigr\}.
	\end{aligned}
\]
By \cite[Theorem~4.1]{KacQua03b}, \smash{$\upr$} is then strongly and freely generated, up to some \opes\ that must be computed, by two even fields of conformal weight $2$ and four odd fields of conformal weights $1$, $1$, $2$ and $2$.

Following the recipe in \cite{KacQua03b}, we associate to each even basis element of \smash{$\psl^{(+2)}$} a fermionic ghost \sva\ $\fgvoa^i$, $i=1,2$, and to each odd basis element a bosonic ghost \va\ $\bgvoa^j$, $j=e,f$.
These are strongly and freely generated by fields satisfying
\[
	b^i(z)c^i(w) \sim \frac{\wun}{z-w} \qquad \text{and} \qquad \beta^j(z)\gamma^j(w) \sim -\frac{\wun}{z-w},
\]
respectively, with all other \opes\ of the generators being regular.
We next form the graded complex
\smash{$\VOA{C} = \upsl \otimes \fgvoa^1 \otimes \fgvoa^2 \otimes \bgvoa^e \otimes \bgvoa^f$},
in which the $b^i$ and $\beta^j$ have grade $-1$, the $c^i$ and $\gamma^j$ have grade $+1$, and all other generators (and vacuum states) have grade $0$.
This becomes a differential graded complex upon taking the differential to be the zero mode of
\[
	D = \big(E^1 + \wun\big) c^1 + \big(E^2 - \wun\big) c^2 - e^{++} \gamma^e - f^{++} \gamma^f.
\]
$D_0$ is odd, has grade $+1$ and is easily checked to square to zero.
As a \sva, \smash{$\upr$} is the grade-$0$ cohomology of $\VOA{C}$ with respect to $D_0$.

To give \smash{$\upr$}, $\kk\ne0$, the structure of a \svoa, we choose an \emt\ $L$ on $\VOA{C}$.
This needs to make $D$ homogeneous of conformal weight $1$, so we take
\[
	L = T + \tfrac{1}{2} \big(\pd H^1 + \pd H^2\big) + \no{\pd b^1 c^1} + \no{\pd b^2 c^2} + \no{\pd \beta^e \gamma^e} + \no{\pd \beta^f \gamma^f}.
\]
It is easy to check that $\VOA{C}$ is now a \svoa\ of central charge $-2$ and that $D_0 L = 0$.
Moreover, $L$ defines an \emt\ on \smash{$\upr$}, as desired, by \cite[Theorem~2.2]{KacQua03}.

The next step is to construct the ``parenthetical building blocks'' of the generators of \smash{$\upr$}.
These are obtained by correcting some of the generators of \smash{$\upsl$} by adding bilinear terms in the ghost fields, see \cite[equation~(2.4)]{KacQua03b}.
We need to compute them for all elements in \smash{$\psl^{(0)}$}, the results for our basis being
\begin{alignat}{3}
		& H^{(1)} = H^1 + 2\no{b^1c^1} + \no{\beta^e\gamma^e} + \no{\beta^f\gamma^f},\qquad&&
		H^{(2)} = H^2 + 2\no{b^2c^2} + \no{\beta^e\gamma^e} + \no{\beta^f\gamma^f}, &\nonumber\\
		&e^{(+-)} = e^{+-} - \no{\beta^ec^2} + \no{b^1\gamma^f},\qquad&&
		e^{(-+)} = e^{-+} - \no{\beta^ec^1} - \no{b^2\gamma^f}, &\nonumber\\
		&f^{(+-)} = f^{+-} - \no{\beta^fc^2} - \no{b^1\gamma^e},\qquad&&
		f^{(-+)} = f^{-+} - \no{\beta^fc^1} + \no{b^2\gamma^e}.&\label{eq:(KW)}
\end{alignat}

The generating fields of conformal weight $1$ in \smash{$\upr$} are then (the images in cohomology of) the linear combinations of the parenthetical building blocks that commute with $F$ \cite[Theorem~4.1]{KacQua03b}.
In this way, we find two odd generators
\smash{$
	\chi^e = e^{(+-)} - e^{(-+)}$} and \smash{$\chi^f = f^{(+-)} - f^{(-+)}$}.
It is easy to check that these are $D_0$-closed.

To determine the weight-$2$ generating fields of \smash{$\upr$}, we must find $D_0$-closed corrections to the elements of \smash{$\psl^{(-2)}$} (all of which commute with $F$).
This time, \cite{KacQua03b} does not give a precise recipe, only guaranteeing that the corrections belong to the subalgebra generated by the parenthetical building blocks \eqref{eq:(KW)}.
Finding the corrections by brute force is quite tractable in this case (we used \textsc{OPEdefs} \cite{ThiOPE91}).
The resulting four generators are
\begin{gather*}
		B^1 = F^1 - \frac{1}{4} \no{H^{(1)}H^{(1)}} - \frac{\kk}{2} \pd H^{(1)} - \no{e^{(-+)}f^{(-+)}}
		\\
\phantom{B^1 =}{}	- \frac{(2\kk-1)(3\kk+1)}{4\kk} \no{\big(e^{(+-)}-e^{(-+)}\big) \big(f^{(+-)}-f^{(-+)}\big)}, \\
		B^2 = F^2 + \frac{1}{4} \no{H^{(2)}H^{(2)}} - \frac{\kk}{2} \pd H^{(2)} - \no{e^{(+-)}f^{(+-)}}
			\\
\phantom{B^2 =}{}+ \frac{(2\kk+1)(3\kk-1)}{4\kk} \no{\big(e^{(+-)}-e^{(-+)}\big) \big(f^{(+-)}-f^{(-+)}\big)}, \\
		\psi^e = e^{--} - \frac{1}{2} \brac[\big]{\no{H^{(1)}e^{(+-)}} - \no{H^{(2)}e^{(-+)}}} - \frac{2\kk-1}{4} \pd e^{(+-)} - \frac{2\kk+1}{4} \pd e^{(-+)}, \\
		\psi^f = f^{--} - \frac{1}{2} \brac[\big]{\no{H^{(1)}f^{(+-)}} - \no{H^{(2)}f^{(-+)}}} - \frac{2\kk-1}{4} \pd f^{(+-)} - \frac{2\kk+1}{4} \pd f^{(-+)},
\end{gather*}
where unique solutions were obtained by requiring in addition that the generators be primary with respect to $L$.
We record that (up to $D_0$-exact terms)
\[
	L = -\frac{1}{\kk} \big(B^1+B^2\big) - \frac{1}{2\kk} \no{\chi^e\chi^f}.
\]

\subsection{Operator product expansions} \label{sec:opes}

Having explicitly computed the generating fields of the universal level-$\kk$ principal W-algebra~$\upr$, we now turn to its \opes.
By \cite[Theorem~4.1 and Remark~4.2]{KacQua03b}, \smash{$\upr$}~admits a \pbw\ basis in the modes of the generating fields.
It follows that~\smash{$\upr$} is strongly and freely generated by $B^1$, $B^2$, $\chi^e$, $\chi^f$, $\psi^e$ and $\psi^f$, subject only to their \opes.
In principle, these relations might only close modulo $D_0$-exact fields.
Happily, this is not the case.

The resulting \opes\ are nevertheless somewhat messy.
To bring out some of the structural features, we first present those of the generating fields of conformal weight~$1$
\[
	\chi^e(z) \chi^e(w) \sim 0, \qquad \chi^e(z) \chi^f(w) \sim \frac{2\kk\wun}{(z-w)^2}, \qquad \chi^f(z) \chi^f(w) \sim 0.
\]
The weight-$1$ fields of \smash{$\upr$} thus generate a subalgebra isomorphic to a pair of symplectic fermions~\cite{KauSym00}.
We may even remove the apparent level dependence by defining renormalised symplectic fermions thusly
$\chi = \frac{1}{\sqrt{\kk}} \chi^e $ and $\bchi = \frac{1}{\sqrt{\kk}} \chi^f$.

To choose a better basis of weight-$2$ generating fields, we require that they transform ``cleanly'' under the action of the symplectic fermion subalgebra.
More precisely, the space of weight-$2$ generators forms a representation of the zero modes of the symplectic fermions \svoa\ $\sfvoa = \uaff{\kk}{\pslone}$ (these are isomorphic for all $\kk\ne0$), so we can look for a basis that respects the structure of this representation.
Since the zero modes of the generators define a~copy of $\pslone$, a Lie superalgebra with relatively few representations, it is in fact rather easy to identify ours.
Indeed, if we define
\[
	S = -\frac{1}{\kk}\big(B^1+B^2\big), \qquad H = \tfrac{1}{2}\big(B^2-B^1\big) - \tfrac{1}{2}S, \qquad
	\psi = \frac{1}{\sqrt{\kk}}\psi^e \qquad \text{and} \qquad \bpsi = \frac{1}{\sqrt{\kk}}\psi^f,
\]
it is straightforward to verify (again, using \textsc{OPEdefs}) that
\begin{gather*}
	\chi_0 H = \psi, \qquad \bchi_0 H = \bpsi, \qquad
	\chi_0 \bpsi = +S, \qquad \bchi_0 \psi = -S \\ \text{and} \qquad
	\chi_0 \psi = \bchi_0 \bpsi = \chi_0 S = \bchi_0 S = 0.
\end{gather*}
The space spanned by the weight-$2$ generators is thus a copy of the unique projective module of~$\pslone$, which is itself isomorphic to \smash{$\envalg{\pslone}$}.

Note that the \emt\ of \smash{$\upr$} is now given by
$L = S - \frac{1}{2} \no{\chi\bchi}$.
Because of this, we may replace $S$ by $L$ in our set of strong generators.
We record the defining \opes\ of \smash{$\upr$} in the following Theorem~\ref{thm:propes}.
\begin{Theorem} \label{thm:propes}
	The principal \qhr\ \smash{$\upr$} of \smash{$\upsl$} is strongly and freely generated by two even elements, $L$ and $H$, and four odd elements, $\chi$, $\bchi$, $\psi$ and $\bpsi$.
	Here, $L$ is a conformal vector of central charge $-2$ and the generators $H$, $\chi$, $\bchi$, $\psi$ and $\bpsi$ are Virasoro primaries with respective conformal weights $2$, $1$, $1$, $2$ and $2$.
	The remaining \opes\ are as follows:
	\begin{subequations} \label{eq:opes}
		\begin{gather}
			\begin{aligned}
				&\chi(z) \chi(w) \sim \bchi(z) \bchi(w) \sim 0, \qquad
				\chi(z) \bpsi(w) \sim +\frac{S(w)}{z-w},
				\\
				&\chi(z) \psi(w) \sim \bchi(z) \bpsi(w) \sim 0, \qquad
				\bchi(z) \psi(w) \sim -\frac{S(w)}{z-w},\qquad
 \chi(z) \bchi(w) \sim \frac{2\wun}{(z-w)^2},
			\end{aligned}
			\\
			\begin{aligned}
				&\chi(z) H(w) \sim -3\brac*{\kk^2-\tfrac{1}{4}} \frac{\chi(w)}{(z-w)^2} + \frac{\psi(w)}{z-w}, \qquad
				\psi(z) \psi(w) \sim \frac{\no{\chi\psi}(w)}{z-w}, \\
				& \bchi(z) H(w) \sim -3\brac*{\kk^2-\tfrac{1}{4}} \frac{\bchi(w)}{(z-w)^2} + \frac{\bpsi(w)}{z-w}, \qquad
				 \bpsi(z) \bpsi(w) \sim \frac{\no{\bchi\bpsi}(w)}{z-w},
			\end{aligned}
			\\
			\begin{aligned}
				H(z) \psi(w) \sim{}& \tfrac{3}{2}\brac*{\kk^2-\tfrac{1}{4}} \sqbrac*{\frac{\chi(w)}{(z-w)^3} + \frac{\frac{1}{2}\pd\chi(w)}{(z-w)^2} + \frac{\frac{1}{6}\pd^2\chi(w)}{z-w}} \\
&+ \frac{\frac{1}{4}\pd\psi(w) - \no{\chi\big(H+\frac{3\kk^2}{2}S\big)}(w)}{z-w}, \\
				H(z) \bpsi(w) \sim{}& \tfrac{3}{2}\brac*{\kk^2-\tfrac{1}{4}} \sqbrac*{\frac{\bchi(w)}{(z-w)^3} + \frac{\frac{1}{2}\pd\bchi(w)}{(z-w)^2} + \frac{\frac{1}{6}\pd^2\bchi(w)}{z-w}}\\
& + \frac{\frac{1}{4}\pd\bpsi(w) - \no{\bchi\big(H+\frac{3\kk^2}{2}S\big)}(w)}{z-w},
			\end{aligned}
			\\
			\begin{aligned}
				H(z) H(w) \sim{} & -3(3\kk^2-1)\brac*{\kk^2-\tfrac{1}{4}} \sqbrac*{\frac{\wun}{(z-w)^4} + \frac{\no{\chi\bchi}(w)}{(z-w)^2} + \frac{\frac{1}{2}\pd\no{\chi\bchi}(w)}{z-w}} \\
				& + \frac{7\kk^2-1}{2} \sqbrac*{\frac{S(w)}{(z-w)^2} + \frac{\frac{1}{2}\pd S(w)}{z-w}},
			\end{aligned}
			\\
			\begin{aligned}
				\psi(z) \bpsi(w) \sim{} & 3\brac*{\kk^2-\tfrac{1}{4}} \sqbrac*{\frac{\wun}{(z-w)^4} + \frac{\no{\chi\bchi}(w)}{(z-w)^2} + \frac{\frac{1}{2}\pd\no{\chi\bchi}(w)}{z-w}}
				\\
				& - \sqbrac*{\frac{2H(w)+\frac{1}{2}S(w)}{(z-w)^2} + \frac{\pd H(w)+\frac{1}{4}\pd S(w)-\frac{1}{2}\no{\chi\bpsi}(w)-\frac{1}{2}\no{\bchi\psi}(w)}{z-w}}.
			\end{aligned}
		\end{gather}
	\end{subequations}
\end{Theorem}

In Theorem~\ref{thm:propes}, we chose to present the defining \opes\ \eqref{eq:opes} so that the \rhss\ involve $S$ rather than $L$, in order to respect the $\pslone$-symmetry.
For completeness, we also give the \opes\ in which the \lhss\ involve~$S$%
\begin{subequations} \label{eq:opesS}
	\begin{gather}
		\chi(z) S(w) \sim \bchi(z) S(w) \sim 0, \qquad
		S(z) S(w) \sim \frac{2S(w)}{(z-w)^2} + \frac{\pd S(w)}{z-w},
		\\
		\begin{aligned}
			S(z) \psi(w) &\sim \frac{2\psi(w)}{(z-w)^2} + \frac{\pd\psi(w) - \frac{1}{2}\no{\chi S}(w)}{z-w}, \\
			S(z) \bpsi(w) &\sim \frac{2\bpsi(w)}{(z-w)^2} + \frac{\pd\bpsi(w) - \frac{1}{2}\no{\bchi S}(w)}{z-w},
		\end{aligned}
		\\
		\begin{aligned}[b]
			S(z) H(w) \sim& -3\brac*{\kk^2-\tfrac{1}{4}}\sqbrac*{\frac{\wun}{(z-w)^4} + \frac{\no{\chi\bchi}(w)}{(z-w)^2} + \frac{\frac{1}{2}\pd\no{\chi\bchi}(w)}{z-w}} \\
			&+ \frac{2H(w) + \frac{1}{2}S(w)}{(z-w)^2} + \frac{\pd H(w) + \frac{1}{2}\no{\chi\bpsi}(w) - \frac{1}{2}\no{\bchi\psi}(w)}{z-w}.
		\end{aligned}
	\end{gather}
\end{subequations}

We make a few additional comments about \smash{$\upr$}:
\begin{itemize}\itemsep=0pt
	\item The \opes\ \eqref{eq:opes} only depend on $\kk^2$, not $\kk$, so that $\upr \cong \upr[-\kk]$.
	This is an obvious consequence of the isomorphism $\upsl \cong \upsl[-\kk]$, noted in Section~\ref{sec:Walg}, and the fact that the principal nilpotent $F^1+F^2 \in \psl$ is invariant under $\omega$.
	\item $S$ is an \emt\ of central charge $0$.
	It commutes with the symplectic fermions subalgebra generated by $\chi$ and $\bchi$.
	It would be interesting to identify the commutant of this subalgebra in \smash{$\upr$}.
	\item When $\kk^2=\frac{1}{4}$ (hence $\kk=\pm\frac{1}{2}$), the level is ``collapsing'' \cite{AdaCon16}, meaning that the simple quotient~$\spr$ of \smash{$\upr$} is strongly generated by just its weight-$1$ fields.
	To see this, note that~$H$,~$S$,~$\psi$ and $\bpsi$ generate an ideal of \smash{$\upr[\pm1/2]$} and that the quotient is the symplectic fermions subalgebra (which is simple): \smash{$\spr[\pm1/2] \cong \sfvoa$}.
	\item \smash{$\upr$} admits a horizontal grading in which $\chi$ and $\psi$ have grade $+1$, $\bchi$ and $\bpsi$ have grade $-1$, and $H$, $S$ and $L$ have grade $0$.
	Its graded character $\big(\chmap^+\big)$ and supercharacter ($\chmap^-$) are thus
	\begin{gather} \label{eq:uprchar}
		\fgch{\upr}{\yy;\qq} = \qq^{1/12} \prod_{i=1}^{\infty} \frac{\big(1\pm\yy\qq^i\big)\big(1\pm\yy^{-1}\qq^i\big)\big(1\pm\yy\qq^{i+1}\big)\big(1\pm\yy^{-1}\qq^{i+1}\big)}{\big(1-\qq^{i+1}\big)^2}.
	\end{gather}
	\item Based on numerical investigations with the analogue of the Shapovalov form, we conjecture that when \smash{$\upr$} is not simple, then the singular vector of minimal conformal weight has grade~$+1$ and multiplicity $1$.
 The conformal weight is $2$ for $\kk=\pm\frac{1}{2}$, $4$ for $\kk=\pm\frac{1}{3}, \pm\frac{3}{2}$, and $6$ for~${\kk=\pm\frac{1}{4}, \pm\frac{2}{3}, \pm\frac{5}{2}}$.
\end{itemize}

The simplicity of \smash{$\upr$} is obviously of interest.
The following result establishes this for all but a few levels.
\begin{Theorem} \label{thm:prsimple}
 The \svoa\ \smash{$\upr$} is simple if $\kk \notin \QQ$ or $\kk \in \ZZ \setminus \{0,\pm1,\pm2\}$.
 It is not simple if $\kk \in \QQ \setminus \ZZ$.
\end{Theorem}
\begin{proof}
 We make use of the Semikhatov realisation $\umin \ira \upr \otimes \lvoa$ of Theorem~\ref{thm:usemi} and (a spectral-flow twist of) an Adamovi\'{c} functor \eqref{eq:adfunct}, see Section~\ref{sec:iqhr} below.
 They connect the (universal) principal and minimal W-algebras of \smash{$\psl$} via inverse \qhr.
 \begin{itemize}\itemsep=0pt
 \item Building on \cite{GorSim06,HoySim08}, Gorelik and Kac prove in \cite[Corollary~1.1\,(iv)]{GorSim23} that $\umin$ is simple if $\kk \notin \QQ$ or $\kk \in \ZZ_{\ge3}$ and that it is not simple if $\kk \in \QQ \setminus \ZZ$.
 \item Assume that \smash{$\upr$} has a nonzero proper submodule $\VOA{I}$.
 Then, applying an Adamovi\'{c} functor shows that $\VOA{I} \otimes \lvoa$ is a nonzero proper $\umin$-submodule of $\upr \otimes \lvoa$.
 The argument of \cite[Lemma~3.14]{AdaWei23} now proves that for any $v \in \VOA{I}$, we have \smash{$v \otimes \ee^{2nc} \in \umin$} for all sufficiently large~${n\in\ZZ}$.
 $\VOA{I} \otimes \lvoa$ is thus a nonzero proper submodule of $\umin$.
 Contrapositively, $\kk \notin \QQ$ or~${\kk \in \ZZ_{\ge3}}$ implies that $\umin$ is simple, hence so are \smash{$\upr$} and $\upr[-\kk]$.
 \item Next, assume that \smash{$\upr$} is simple \emph{and} that \smash{$\kk \notin \ZZ_{\le0}$}.
 Then, \smash{$\upr \otimes \lvoa$} is (up to a spectral flow) an almost-irreducible $\umin$-module containing $\umin$.
 But, this forces $\umin$ to be simple because the restriction on $\kk$ rules out a singular vector of the form $\wun \otimes \ee^{2nc}$, $n \in \ZZ_{>0}$.
 Contrapositively, $\kk \in \QQ \setminus \ZZ$ implies that $\umin$ is not simple, hence that \smash{$\upr$} is not simple either.\hfill $\qed$
 \end{itemize}\renewcommand{\qed}{}
\end{proof}

 We thank Drazen Adamovi{\'c} for pointing out to us the work \cite{GorSim23} leading to the proof of Theorem~\ref{thm:prsimple}.
The simplicity of \smash{$\upr$} is not settled for $\kk=\pm1,\pm2$ ($\kk=0$ is critical), but we expect a positive answer in these cases.
We also thank an anonymous reviewer for sketching an alternative, and extremely interesting, approach to proving this result when $\kk \in \QQ \setminus \ZZ$.

\section[Irreducible upr-modules]{Irreducible \smash{$\boldsymbol{\upr}$}-modules} \label{sec:reps}

We next study the representation theory of \smash{$\upr$}.
Being a principal W-algebra, we expect that its irreducible modules are all ordinary, meaning that the $L_0$-eigenspaces are all \fdim.
This is proven under the additional hypothesis that the irreducible is lower bounded.
We also conjecture that every irreducible \smash{$\upr$}-module is lower bounded.

\subsection{The principal Zhu algebra} \label{sec:zhu}

The Zhu algebra $\zhu{\VOA{V}}$ of a \svoa\ $\VOA{V}$ is the unital associative superalgebra of zero modes of the fields, restricted to act on states that are annihilated by all positive modes.
Introduced formally in \cite{KacVer93,ZhuMod96}, $\zhu{\VOA{V}}$ may be realised as a quotient of $\VOA{V}$ with an associative product given by
\begin{equation} \label{eq:zhudef}
	\zh{A} \zh{B} = \sum_{j=0}^{\Delta_A} \binom{\Delta_A}{j} \zh{(A_{-\Delta_A+j}B)}, \qquad A,B \in \VOA{V}.
\end{equation}
Here, $\zh{A}$ denotes the image of $A$ in $\zhu{\VOA{V}}$ and $\Delta_A$ is the conformal weight of $A$ (which is assumed to be homogeneous).
It inherits its parity from that of $A$.
The quotient operation will not be important for us here, except to note that it includes the identification
\begin{equation} \label{eq:zhuid}
	\zh{(\pd A)} = -\Delta_A \zh{A}, \qquad A \in \VOA{V}.
\end{equation}

It is now easy to use \eqref{eq:zhudef} and \eqref{eq:zhuid} to deduce relations in \smash{$\zhu{\upr}$} involving the images of the strong generators
\begin{gather}
			(\zh{\chi})^2 = (\zh{\bchi})^2 = 0,\qquad
			\acomm[\big]{\zh{\chi}}{\zh{S}} = \acomm[\big]{\zh{\bchi}}{\zh{S}} = 0,\qquad
			\acomm[\big]{\zh{\chi}}{\zh{H}} = \zh{\psi},\qquad
			\acomm[\big]{\zh{\bchi}}{\zh{H}} = \zh{\bpsi},\nonumber
			\\
			\acomm[\big]{\zh{\chi}}{\zh{\bchi}} = 0,\qquad
			\acomm[\big]{\zh{\chi}}{\zh{\psi}} = \acomm[\big]{\zh{\bchi}}{\zh{\bpsi}} = 0,\qquad
			\acomm[\big]{\zh{\chi}}{\zh{\bpsi}} = +\zh{S},\qquad
			\acomm[\big]{\zh{\bchi}}{\zh{\psi}} = -\zh{S},\nonumber
		\\
		(\zh{\psi})^2 = \tfrac{1}{2}\zh{\chi}\zh{\psi}, \qquad
		\big(\zh{\bpsi}\big)^2 = \tfrac{1}{2}\zh{\bchi}\zh{\bpsi}, \qquad
		\comm[\big]{\zh{\psi}}{\zh{S}} = \tfrac{1}{2}\zh{\chi}\zh{S}, \qquad
		\comm[\big]{\zh{\bpsi}}{\zh{S}} = \tfrac{1}{2}\zh{\bchi}\zh{S},\label{eq:zhurels}
		\\
			\acomm[\big]{\zh{\psi}}{\zh{\bpsi}} = \tfrac{1}{2}\brac[\big]{\zh{\chi}\zh{\bpsi} + \zh{\bchi}\zh{\psi}}, \qquad
			\comm[\big]{\zh{\psi}}{\zh{H}} = \zh{\chi}\brac[\bigg]{\zh{H}+\frac{3\kk^2}{2}\zh{L}+\tfrac{1}{4}\big(\kk^2-\tfrac{1}{4}\big)\zh{\wun}} - \tfrac{1}{2}\zh{\psi},\nonumber
			\\
			\comm[\big]{\zh{S}}{\zh{H}} = \tfrac{1}{2}\brac[\big]{\zh{\chi}\zh{\bpsi} + \zh{\psi}\zh{\bchi}}, \qquad
			\comm[\big]{\zh{\bpsi}}{\zh{H}} = \zh{\bchi}\brac[\bigg]{\zh{H}+\frac{3\kk^2}{2}\zh{L}+\tfrac{1}{4}\bigl(\kk^2-\tfrac{1}{4}\bigr)\zh{\wun}} - \tfrac{1}{2}\zh{\bpsi}.\nonumber
\end{gather}
We shall not do so here, but one can prove that the images of the strong generators, along with the relations \eqref{eq:zhurels}, comprise a presentation of \smash{$\zhu{\upr}$}.
Instead, we record some useful observations:
\begin{itemize}\itemsep=0pt
	\item As usual, $\zh{L} = \zh{S} + \frac{1}{2} \zh{\bchi}\zh{\chi}$ is central.
	\item The Zhu-images of the odd weight-$2$ fields are not nilpotent of order $2$, as one might expect, but of order $3$:
$(\zh{\psi})^3 = \frac{1}{2} \zh{\chi}(\zh{\psi})^2 = \frac{1}{4} (\zh{\chi})^2\zh{\psi} = 0$.
	\item The horizontal grading of \smash{$\upr$} descends to a grading of \smash{$\zhu{\upr}$} in which $\zh{\chi}$ and $\zh{\psi}$ have grade $+1$, $\zh{\bchi}$ and $\zh{\bpsi}$ have grade $-1$, and $\zh{H}$, $\zh{L}$ and $\zh{S}$ have grade $0$.
	\item A straightforward consequence of the relations \eqref{eq:zhurels} is that monomials in the generators of~$\zhu{\upr}$ may be consistently ordered so that the odd elements with bars (no bars) are placed to the left (right) of the even elements.
	For example, the juxtaposition $\zh{\psi} \zh{H}$ may be so reordered by writing
	\begin{align*}
			\zh{\psi} \zh{H}
			&= \zh{H}\zh{\psi} + \zh{\chi}\brac[\bigg]{\zh{H}+\frac{3\kk^2}{2}\zh{L}+\tfrac{1}{4}\bigl(\kk^2-\tfrac{1}{4}\bigr)\zh{\wun}} - \tfrac{1}{2}\zh{\psi} \\
			&= \bigl(\zh{H} + \tfrac{1}{2}\zh{\wun}\bigr) \zh{\psi} + \brac[\bigg]{\zh{H}+\frac{3\kk^2}{2}\zh{L}+\tfrac{1}{4}\bigl(\kk^2-\tfrac{1}{4}\bigr)\zh{\wun}} \zh{\chi}.
	\end{align*}
	This \pbw-style ordering suggests that \smash{$\zhu{\upr}$} admits a well behaved \hw\ theory.
\end{itemize}

\subsection{Modules of the Zhu algebra} \label{sec:zhuwt}

Let us define a weight vector in a \smash{$\zhu{\upr}$}-module to be a simultaneous eigenvector of $\zh{H}$ and~$\zh{L}$.
As usual, weight vectors always exist in (nonzero) \fdim\ modules.
We will call a weight vector a \hwv\ if it is annihilated by both $\zh{\chi}$ and $\zh{\psi}$.
Obviously, a~\hwv\ is also an eigenvector of $\zh{S}$ (and the eigenvalue matches that of~$\zh{L}$).\looseness=1

\begin{Lemma} \label{lem:zhuwtmod}
	If a \smash{$\zhu{\upr}$}-module has a weight vector, then it has a \hwv.
\end{Lemma}
\begin{proof}
	Let $v$ denote the hypothesised weight vector, so that $\zh{L} v = \Delta v$ and $\zh{H} v = hv$, for some~${\Delta,h \in \CC}$.
	Then,
	\begin{itemize}\itemsep=0pt
		\item If $\zh{\chi} v = \zh{\psi} v = 0$, then $v$ is already \hw.
		\item If $\zh{\psi} v = 0$ but $\zh{\chi} v \ne 0$, then $\zh{\chi} v$ is weight because $\zh{L}$ is central and $\zh{H}\zh{\chi} v = \zh{\chi}\zh{H} v - \zh{\psi} v = h \zh{\chi} v$.
		But, $\zh{\chi}\zh{\chi} v = 0$ and $\zh{\psi}\zh{\chi} v = -\zh{\chi}\zh{\psi} v = 0$, so $\zh{\chi} v$ is \hw.
		\item If $\zh{\chi} v = 0$ but $\zh{\psi} v \ne 0$, then $\zh{\psi} v$ is weight because
		\[
			\zh{H}\zh{\psi} v = \zh{\psi}\zh{H} v + \tfrac{1}{2} \zh{\psi} v - \zh{\chi}\zh{H} v - \frac{3\kk^2}{2} \zh{\chi}\zh{L} v - \tfrac{1}{4}\bigl(\kk^2-\tfrac{1}{4}\bigl)\zh{\chi} v = \bigl(h+\tfrac{1}{2}\bigr) v.
		\]
		But, $\zh{\chi}\zh{\psi} v = -\zh{\psi}\zh{\chi} v = 0$ and $\zh{\psi}\zh{\psi} v = \frac{1}{2} \zh{\chi}\zh{\psi} v = 0$, so $\zh{\psi} v$ is \hw.
		\item If $\zh{\chi} v \ne 0$ and $\zh{\psi} v \ne 0$ but $\zh{\chi}\zh{\psi} v = 0$, then
		\[
			\zh{H}\zh{\chi} v = h \zh{\chi} v - \zh{\psi} v, \qquad
			\zh{H}\zh{\psi} v = \bigl(h+\tfrac{1}{2}\bigr) \zh{\psi} v - \brac[\bigg]{h + \frac{3\kk^2}{2}\Delta + \tfrac{1}{4}\bigl(\kk^2-\tfrac{1}{4}\bigr)} \zh{\chi} v.
		\]
		Some linear combination of $\zh{\chi} v$ and $\zh{\psi} v$ is therefore a weight vector.
		Since $\zh{\chi} \zh{\psi} v = 0$, this linear combination is annihilated by both $\zh{\chi}$ and $\zh{\psi}$, hence it is a \hwv.
		\item Finally, if $\zh{\chi}\zh{\psi} v \ne 0$ (so $\zh{\chi} v \ne 0$ and $\zh{\psi} v \ne 0$), then $\zh{H}\zh{\chi} \zh{\psi} v = h \zh{\chi} \zh{\psi} v$ and $\zh{\chi} \zh{\chi} \zh{\psi} v = \zh{\psi} \zh{\chi} \zh{\psi} v = 0$, so $\zh{\chi} \zh{\psi} v$ is the desired \hwv. \hfill $\qed$
	\end{itemize} \renewcommand{\qed}{}
\end{proof}

Inspection of \eqref{eq:zhurels} demonstrates that $\zh{H}$, $\zh{L}$, $\zh{\chi}$ and $\zh{\psi}$ generate a unital subalgebra of~$\zhu{\upr}$.
We may therefore construct Verma modules in the usual way and realise all \hwms, these being modules generated by a single \hwv, as quotients.
Let \smash{$\zhver{h,\Delta}$} denote the Verma module generated by a \hwv\ of $\zh{H}$\nobreakdash-eigenvalue $h$ and $\zh{L}$\nobreakdash-eigenvalue $\Delta$.

\begin{Proposition} \label{prop:vermas}
\
	\begin{itemize}\itemsep=0pt
		\item \smash{$\zhver{h,\Delta}$} is $4$-dimensional, for all $h,\Delta \in \CC$.
		\item \smash{$\zhver{h,\Delta}$} is irreducible unless $\Delta=0$.
		\item The irreducible quotient of $\zhver{h,0}$ is $1$-dimensional, for all $h\in\CC$.
	\end{itemize}
\end{Proposition}
\begin{proof}
	Let $v$ be the generating \hwv\ of \smash{$\zhver{h,\Delta}$}.
	Because we can order monomials so that the ``negative'' generators $\zh{\bchi}$ and $\zh{\bpsi}$ act last, \smash{$\zhver{h,\Delta}$} is spanned by the monomials obtained by acting on $v$ with these negative generators.
	It is easy to check that these monomials are $v$, $\zh{\bchi} v$, $\zh{\bpsi} v$ and $\zh{\bchi} \zh{\bpsi} v$, hence that \smash{$\zhver{h,\Delta}$} is $4$-dimensional.

	To investigate irreducibility, consider first the submodule of \smash{$\zhver{h,\Delta}$} generated by the linear combination $w = a \zh{\bchi} v + b \zh{\bpsi} v$, where $a,b \in \CC$ are not both $0$.
	Using \eqref{eq:zhurels}, we find that~${\zh{\chi} w = \Delta bv}$ and $\zh{\psi} w = \Delta \bigl(-a+\frac{1}{2}b\bigr) v$.
	When $\Delta \ne 0$, it follows that this submodule contains $v$, for any~${(a,b) \ne (0,0)}$, and is therefore \smash{$\zhver{h,\Delta}$} itself.
	Similarly, $\Delta \ne 0$ implies that the submodule generated by $\zh{\bchi}\zh{\bpsi} v$ contains $\zh{\chi}\zh{\bchi}\zh{\bpsi} v = -\Delta \zh{\bchi} v$ and thus $v$, hence is also \smash{$\zhver{h,\Delta}$}.
	This proves the irreducibility of Verma modules when $\Delta \ne 0$.

	When $\Delta = 0$, however, these calculations show that every nonzero linear combination of~$\zh{\bchi} v$ and $\zh{\bpsi} v$ is a \hwv, hence generates a proper submodule.
	The sum of these proper submodules obviously contains $\zh{\bchi} \zh{\bpsi} v$, so we conclude that the irreducible quotient of~$\zhver{h,0}$ is spanned by (the image of) $v$.
\end{proof}

An irreducible \smash{$\zhu{\upr}$}-module with a weight vector is thus a \hwm\ and so is either $4$- or $1$-dimensional.
There are nevertheless reducible but indecomposable \smash{$\zhu{\upr}$}-modules, also possessing weight vectors, of dimensions greater than $4$.
We shall not describe these here, but note instead that the action of \smash{$\zhu{\upr}$} on Verma modules is not as nice as we might have hoped.
In particular, the action of $\zh{S}$ on \smash{$\zhver{h,\Delta}$}, $\Delta\ne0$, is nondiagonalisable
\begin{alignat*}{3}
& \zh{S} v = \Delta v, \qquad&&
		\zh{S}\zh{\bchi}\zh{\bpsi} v = \Delta \zh{\bchi}\zh{\bpsi} v,&
\\
&\zh{S}\zh{\bchi} v = \Delta \zh{\bchi} v,\qquad&&
		\zh{S}\zh{\bpsi} v = -\frac{1}{2} \Delta \zh{\bchi} v + \Delta \zh{\bpsi} v.&
\end{alignat*}
The action of $\zh{H}$ is even more opaque
\begin{alignat*}{3}
		&\zh{H} v = hv,\qquad&&
		\zh{H}\zh{\bchi}\zh{\bpsi} v = h \zh{\bchi}\zh{\bpsi} v,&
\\
		&\zh{H}\zh{\bchi} v = h \zh{\bchi} v - \zh{\bpsi} v, \qquad&&
		\zh{H}\zh{\bpsi} v = -\brac[\bigg]{h + \frac{3\kk^2}{2}\Delta + \tfrac{1}{4}\bigl(\kk^2-\tfrac{1}{4}\bigr)} \zh{\bchi} v + \bigl(h+\tfrac{1}{2}\bigr) \zh{\bpsi} v.&
\end{alignat*}
The eigenvalues of $\zh{H}$ on the subspace spanned by $\zh{\bchi} v$ and $\zh{\bpsi} v$ are thus
\[
h+\tfrac{1}{4} \pm \sqrt{h+\tfrac{1}{4}(6\Delta+1)\kk^2}.
\]
This action is therefore also nondiagonalisable when $h = -\frac{1}{4}(6\Delta+1)\kk^2$.

This classifies irreducible \smash{$\zhu{\upr}$}-modules with a weight vector.
Could there exist an irreducible \smash{$\zhu{\upr}$}-module without a weight vector?
To tackle this question, recall that because \smash{$\zhu{\upr}$} has countable dimension, its irreducible modules do too.
Dixmier's lemma therefore shows that $\zh{L}$ always acts as a multiple of the identity on any irreducible.

Choose $v\ne0$ in an irreducible \smash{$\zhu{\upr}$}-module.
By acting with $\zh{\chi}$ and $\zh{\psi}$, we may assume without loss of generality that $v$ is annihilated by both $\zh{\chi}$ and $\zh{\psi}$.
If $v$ is a weight vector, then our irreducible is a \hwm\ by Lemma~\ref{lem:zhuwtmod}.
So suppose that $v$ is not weight.
As~$\zh{H} v$ is then nonzero, the submodule it generates must still contain $v$, by irreducibility.
Because we can order monomials in \smash{$\zhu{\upr}$}, this submodule is spanned by monomials of the form~\smash{$\big(\zh{\bpsi}\big)^{\ell} (\zh{\bchi})^m \zh{H}^n v$}, with $n\ge1$.
However, our irreducible admits a grading inherited from that of~\smash{$\zhu{\upr}$} by setting the grade of $v$ to $0$.
To obtain $v$, we may thus restrict to monomials with~${\ell=m=0}$, concluding that $v = p(\zh{H}) v$ for some polynomial $p$ with zero constant term.
But now, $\zh{H}$ preserves the \fdim\ subspace $\{\zh{H}^n v \mid n=0,1,\dots,\deg(p)-1\}$ and so this subspace contains a weight vector.
We therefore conclude that every irreducible has a~weight vector.
\begin{Proposition} \label{prop:irreps_are_hw}
	Every irreducible \smash{$\zhu{\upr}$}-module is highest weight.
\end{Proposition}

\subsection[Irreducible lower-bounded upr- and spr-modules]{Irreducible lower-bounded \smash{$\boldsymbol{\upr}$}- and $\boldsymbol{\spr}$-modules} \label{sec:pslirreps}

Recall that a finitely generated module over a \svoa\ is lower bounded if its conformal weights have a minimal real part.
The subspace realising this minimum is called the top space, and it is naturally a module for the Zhu algebra.
We define a weight vector in a~\smash{$\upr$}-module to be an eigenvector of $H_0$ that is also a generalised eigenvector of $L_0$, its weight then being the pair $(h,\Delta)$ of eigenvalues of $H_0$ and $L_0$ (respectively).
A weight vector is a~\hwv\ if it is also annihilated by $\chi_0$, $\psi_0$ and all the modes of \smash{$\apsl$} with positive indices.
The previous Proposition~\ref{prop:irreps_are_hw} now has the following consequence for the representation theory of \smash{$\upr$}.
\begin{Theorem} \label{thm:irreps_are_hw}
	Every irreducible lower-bounded \smash{$\upr$}-module is highest weight, so is determined up to isomorphism by its highest weight $(h,\Delta)$.
	Moreover, the top space of the irreducible is $1$-dimensional, if $\Delta = 0$, and is otherwise $4$-dimensional.
\end{Theorem}

The obvious next step is to try to classify the irreducible lower-bounded $\spr$-modules in terms of their highest weights.
For general $\kk$, this is clearly beyond the scope of the paper.
However, it is interesting to investigate whether the number of irreducibles is finite.
When \smash{$\kk=\pm\frac{1}{2}$}, there is only one irreducible as \smash{$\spr[\pm1/2] \cong \sfvoa$}.
We record the results of our investigations for the two next-most-accessible levels: $\kk=\pm\frac{1}{3}$ and $\kk=\pm\frac{3}{2}$:
\begin{itemize}\itemsep=0pt
 \item There is a singular vector in \smash{$\upr[\pm1/3]$} of charge $+1$ and conformal weight $4$ (see Section~\ref{sec:opes}).
 It thus has two linearly independent descendants of charge $0$ obtained by acting with $\bchi_0$ and~$\bpsi_0$, respectively.
 Their zero modes act as scalars on a \hwv\ of weight $(h,\Delta)$ and this action must vanish if the vector is to generate a $\spr[\pm1/3]$-module.
 Explicit calculation now gives two equations in two unknowns and a finite number of solutions
 \begin{gather*}
 \Delta (24h - 5\Delta + 5) = 0 \qquad \text{and} \qquad
 576h^2 - 48h - 120h\Delta - 49\Delta - 63\Delta^2 = 0 \\
 \qquad\Ra \quad (h,\Delta) = (0,0),\bigl(\tfrac{1}{12},0\bigr),\bigl(-\tfrac{5}{9},-\tfrac{5}{3}\bigr),\bigl(-\tfrac{5}{36},\tfrac{1}{3}\bigr).
 \end{gather*}
 $\spr[\pm1/3]$ therefore has finitely many lower-bounded irreducibles.
 We conjecture that it is a~log-rational \svoa.
 \item A similar singular vector in \smash{$\upr[\pm3/2]$} likewise leads to two equations in two unknowns.
 Our expectation that we again have finitely many solutions is however dashed as explicit calculation gives
 \begin{gather*}
 \Delta (4h + 5\Delta + 2) = 0 \qquad \text{and} \qquad
 h (4h + 5\Delta + 2) = 0 \\
 \qquad \Ra \quad (h,\Delta) = (0,0), \bigl(h,-\tfrac{2}{5}(2h+1)\bigr), \qquad h \in \CC.
 \end{gather*}
 While this does not prove that $\spr[\pm3/2]$ has infinitely many irreducibles, we conjecture that it does and so is not log-rational.
\end{itemize}

\section[Inverse reduction and the N=4 superconformal algebra]{Inverse reduction and the $\boldsymbol{N=4}$ superconformal algebra} \label{sec:InvRed}

We conclude by explaining how the paradigm of inverse \qhr\ can be used to transfer representation-theoretic information between W-algebras.
For this, we explicitly construct a so-called Semikhatov realisation of the small $N=4$ superconformal \svoa\ $\umin$ from \smash{$\upr$} and sketch how to determine when it descends to $\smin$ and $\spr$.
As the representation theory of \smash{$\upr$} and $\spr$ is still mysterious for general levels $\kk$, we then restrict to $\kk=\pm\frac{1}{2}$.
The corresponding $N=4$ central charges are $-3$ and $-9$, both of which have been previously considered in the literature.
In particular, inverse reduction allows us to easily recover the $\cc=-9$ results of \cite{AdaRea14,AdaVer20}.
We take the opportunity to complement these known results with similar ones for $\cc=-3$ and, for both central charges, detailed structures for certain nonsemisimple \smash{$\sminpmhalf$}-modules, some of which are logarithmic.

\subsection{Inverse \qhr} \label{sec:iqhr}

As previously mentioned, the \qhr\ $\umin$ of \smash{$\upsl$} with respect to a minimal nilpotent element is isomorphic to the small $N=4$ superconformal \svoa.
In the framework of inverse \qhr\ \cite{AdaRea17,SemInv94}, we expect a~Semikhatov realisation of $\umin$ in terms of \smash{$\upr$} and some free field algebra $\VOA{F}$,
$\umin \hookrightarrow \upr \otimes \VOA{F}$.
This embedding should moreover descend to an embedding involving the simple quotients $\smin$ and $\spr$, at least for most $\kk$, and thereby provide nontrivial insights into the representation theory of $\umin$ and $\smin$.

To construct such a Semikhatov realisation, we will utilise an explicit presentation of $\umin$, following \cite{KacQua03,KacQua03b}.
This has the following features:
\begin{itemize}\itemsep=0pt
	\item A strong generating set consists of even elements $J^+$, $J^0$, $J^-$ and $T$, as well as odd elements $G^+$, $G^-$, $\bG^+$ and $\bG^-$.
	\item $T$ is a conformal vector for $\umin$ of central charge $-6(\kk+1)$ and the corresponding conformal weights of the strong generators labelled by $J$ are $1$, while those of the strong generators labelled by $G$ are $\frac{3}{2}$.
	\item All the strongly generating fields except $T(z)$ are primary with respect to this conformal structure.
	\item The zero modes of $J^+$, $J^0$ and $J^-$ form an $\sltwo$-triple and the vertex subalgebra they generate is isomorphic to the affine vertex algebra $\usl[-\kk-1]$.
	\item $G^+$ and $G^-$ span a fundamental module of $\sltwo$, as do $\bG^+$ and $\bG^-$.
	Here, the labels $+$ and $-$ correspond to the eigenvalue of $J^0$.
	We normalise these $\sltwo$-actions so that $J^{\pm}_0 G^{\mp} = G^{\pm}$ and~${J^{\pm}_0 \bG^{\mp} = \bG^{\pm}}$.
\end{itemize}
Given this information, it only remains to specify the \opes\ involving $G^{\pm}$ and $\bG^{\pm}$, the nonregular of which are given by
\begin{gather*}
 G^{\pm} (z) \bG^{\pm} (w) \sim \pm \frac{2 J^{\pm} (w)}{(z-w)^2} \pm \frac{\pd J^{\pm} (w)}{z-w}, \\
 G^{\pm} (z) \bG^{\mp} (w) \sim \pm \frac{2(\kk+1) \wun}{(z-w)^3} - \frac{J^0 (w)}{(z-w)^2} - \frac{\frac{1}{2} \pd J^0(w) \pm T(w)}{z-w}.
\end{gather*}

We remark that if one requires that the action of the even generating fields is invariant under a $2\pi$-rotation about the origin, while that of the odd fields need only be invariant under a $4\pi$\nobreakdash-rotation, then the representation theories of $\umin$ and $\smin$ decompose into two sectors: the \ns\ sector in which all odd fields are invariant and the Ramond sector in which all odd fields are negated under a $2\pi$-rotation.
These two sectors are moreover equivalent as categories because there exist spectral flow maps $\nsfsymb{\ell}$, parametrised for $\ell \in \frac{1}{2}\ZZ$ by $\ell J^0$ (as in \cite{LiPhy97}).
At the level of the modes of the generating fields, $\nsfsymb{\ell}$ adds $\mp\ell$ to the indices of the $G^{\pm}$- and $\bG^{\pm}$-modes, while it adds $\mp2\ell$ to those of the $J^{\pm}$-modes.
Otherwise, we have
\begin{equation} \label{eq:N=4sf}
	\nsf{\ell}[\big]{J^0_n} = J^0_n + 2(\kk+1)\ell\delta_{n,0}\wun \qquad \text{and} \qquad
	\nsf{\ell}{T_n} = T_n - \ell J^0_n- (\kk+1)\ell^2\delta_{n,0}\wun.
\end{equation}
There is also the conjugation automorphism $\nconjsymb$ of the mode algebra of $\umin$.
It preserves the mode indices and $T$, while acting as the Weyl reflection of $\sltwo$ on $J^0$ and $J^{\pm}$.
The action on the odd generators is given by
$
	\nconj[\big]{G^{\pm}} = G^{\mp} $ and $\nconj[\big]{\bG^{\pm}} = -\bG^{\mp}$.

Let $\VOA{F}$ be the \fms\ bosonisation $\lvoa$ of the bosonic ghosts \va, also known as the half-lattice vertex algebra.
This has strong generators $c$, $d$ and $\ee^{mc}$, for $m \in \ZZ$, whose nonregular \opes\ are \cite{FriCon86}
\[
	c(z) d(w) \sim \frac{2 \wun}{(z-w)^2}, \qquad
	d(z) \ee^{mc}(w) \sim \frac{2m \ee^{mc}(w)}{z-w}.
\]
We also introduce the following convenient alternative basis for the Heisenberg elements in $\lvoa$
$a = +\frac{1}{2} (\kk+1)c + \frac{1}{4} d$, $
	b = -\frac{1}{2} (\kk+1)c + \frac{1}{4} d$.
It is easy to check that $a(z) b(w)$ is regular.
Finally, we equip $\lvoa$ with the conformal vector
$t = \frac{1}{2} \no{cd} - \pd a$.
It corresponds to central charge $-2(3\kk+2)$.
The conformal weights of $c$, $d$ and~$\ee^{mc}$ are $1$, $1$ and~$\frac{1}{2}m$, respectively.

Given this setup, it is easy to verify (for example, using \textsc{OPEdefs}) the following result.
\begin{Theorem} \label{thm:usemi}
	For $\kk \neq 0$, there is a homomorphism $\Phi^{\kk} \colon \umin \to \upr \otimes \lvoa$ of \voas, determined by
	\begin{gather}
			\Phi^{\kk}\big(J^+\big) = \ee^{2c}, \qquad
			\Phi^{\kk}\big(J^0\big) = 2b, \qquad
			\Phi^{\kk}(T) = L + t, \qquad
			\Phi^{\kk}\big(G^+\big) = \no{\chi \ee^c}, \nonumber \\
			\Phi^{\kk}\big(\bG^+\big) = \no{\bchi \ee^c},\qquad
			\Phi^{\kk}\big(G^-\big) = -\no{\brac[\big]{\psi + \chi a - \tfrac{1}{2}\bigl(\kk+\tfrac{1}{2}\bigr)\pd \chi} \ee^{-c}}, \nonumber \\
			\Phi^{\kk}(\bG^-) = -\no{\brac[\big]{\bpsi + \bchi a - \tfrac{1}{2}\bigl(\kk+\tfrac{1}{2}\bigr)\pd \bchi} \ee^{-c}},\nonumber \\
			\Phi^{\kk}(J^-) = \no{\brac[\big]{H - \tfrac{1}{2}(\kk-1) S - aa - \tfrac{1}{2}\bigl(\kk+\tfrac{1}{2}\bigr)(3\kk-1) \chi \bchi + \kk \pd a} \ee^{-2 c}}.\label{eq:usemi}
		\end{gather}
\end{Theorem}

To convert this into a Semikhatov realisation of $\umin$, we need to demonstrate that $\Phi^{\kk}$ is injective.
This follows from a combinatorial argument based on the ``triangularity'' of $\Phi^{\kk}$, see~\cite{AdaRea20} for an example of such an argument.
Here, triangularity refers to the fact that the strong generators of $\umin$ and $\upr \otimes \lvoa$ may be ordered so that the image of the $n$-th generator of $\umin$ is expressed in terms of the $m$-th generators, with $m \le n$.
(For $\lvoa$, we only need to consider $d$ and $\ee^c$ for this purpose, see \cite{AdaRea20}.)
We omit the details.

It remains to obtain a Semikhatov realisation of $\smin$.
As in \cite{AdaRea20}, this follows from a remarkable property of inverse \qhr.
First, recall \cite{BerRep01} that $\lvoa$ admits a family of irreducible lower-bounded modules
\smash{$
	\lmod{\lambda} = \Pi \ee^{-a+\lambda c}$}, $ [\lambda] \in \CC/\ZZ$.
The top space of \smash{$\lmod{\lambda}$} is spanned by the $\ee^{-a+\mu c}$ with $\mu \in [\lambda]$ and it has conformal weight $-\frac{1}{4}(\kk+1)$.
Note that because $\ee^c$ has conformal weight $\frac{1}{2}$, the \smash{$\lmod{\lambda}$} are actually $\ZZ_2$-twisted $\lvoa$-modules.

With these (twisted) $\lvoa$-modules, we introduce the Adamovi\'{c} functors
\begin{equation} \label{eq:adfunct}
	\adfunc{\lambda} = \res^{\upr \otimes \lvoa}_{\umin} (\blank \otimes \lmod{\lambda}), \qquad [\lambda] \in \CC/\ZZ.
\end{equation}
The remarkable property is now that if $\Mod{M}$ is a simple \smash{$\upr$}-module, then $\adfunc{\lambda}(\Mod{M})$ is almost irreducible.
This means that it is lower bounded, generated by its top space, and is such that every nonzero submodule has nonzero intersection with the top space.
We will also omit the details that establish this property, referring again to \cite{AdaRea20} for a similar example.

This property of Adamovi\'{c} functors has many important applications.
One such is the precise condition under which $\Phi^{\kk}$ descends to a Semikhatov realisation
\begin{equation} \label{eq:ssemi}
	\Phi_{\kk} \colon\ \smin \ira \spr \otimes \lvoa
\end{equation}
of the simple quotients.
The idea is as follows \cite{AdaRea20}.
Since $\spr$ is irreducible, $\adfunc{\lambda}\big(\spr\big)$ is an almost-irreducible \smash{$\umin$}-module, hence any of its nonzero submodules intersect its top space nontrivially.
But, for a suitable choice of $[\lambda]$, this module is a spectral flow of $\spr \otimes \lvoa$.
It follows that the image of the maximal submodule of $\umin$ is either $0$ or it has nonzero intersection with the ``spectrally flown top space'' $\spn\{(J^+)^n \ket{0} \mid n \in \ZZ\}$ of $\spr \otimes \lvoa$.
As $J^+ \in \usl[-\kk-1]$, the latter is impossible unless $\kk \in \ZZ_{\le-1}$.
\begin{Theorem} \label{thm:ssemi}
	For $\kk \notin \ZZ_{\le0}$, the embedding of \eqref{eq:usemi} induces the simple Semikhatov realisation~\eqref{eq:ssemi}.
\end{Theorem}

\subsection[Constructing weight sminpmhalf-modules]{Constructing weight $\boldsymbol{\sminpmhalf}$-modules} \label{sec:k=+-1/2}

Recall that the simple principal W-algebra $\spr$ collapses to the symplectic fermions \svoa\ $\sfvoa$ when the level is $\kk=\pm\frac{1}{2}$.
The corresponding simple minimal W-algebras $\smin[1/2]$ and $\smin[-1/2]$ have central charges $-9$ and $-3$, respectively.
For convenience, we specialise the Semikhatov realisation of Theorems~\ref{thm:usemi} and~\ref{thm:ssemi} to these levels by setting $H$, $S$, $\psi$ and $\bpsi$ to $0$.
\begin{Corollary} \label{cor:pmhalf}
	For $\kk = \pm\frac{1}{2}$, there is a Semikhatov realisation $\Phi_{\kk} \colon \smin \ira \sfvoa \otimes \lvoa$ given by
	\begin{alignat}{3}
 & \Phi_{\kk}\big(J^+\big) = \ee^{2c}, \qquad&&
 \Phi_{\kk}(J^-) = -\no{\brac*{aa - \kk \pd a} \ee^{-2c}} - \tfrac{1}{2}\bigl(\kk+\tfrac{1}{2}\bigr)(3\kk-1) \no{\chi \bchi} \ee^{-2c},&\nonumber \\
& \Phi_{\kk}\big(J^0\big) = 2b, \qquad&&
 \Phi_{\kk}(T) = \tfrac{1}{2} \no{cd} - \pd a - \tfrac{1}{2} \no{\chi\bchi},& \nonumber\\
& \Phi_{\kk}\big(G^+\big) = \chi \ee^c, \qquad&&
 \Phi_{\kk}(G^-) = \tfrac{1}{2}\bigl(\kk+\tfrac{1}{2}\bigr)\pd \chi \ee^{-c} - \chi \no{a \ee^{-c}},& \nonumber\\
& \Phi_{\kk}\big(\bG^+\big) = \bchi \ee^c, \qquad&&
 \Phi_{\kk}(\bG^-) = \tfrac{1}{2}\bigl(\kk+\tfrac{1}{2}\bigr)\pd \bchi \ee^{-c} - \bchi \no{a \ee^{-c}}.&\label{eq:simplequotemb}
	\end{alignat}
\end{Corollary}

These free-field realisations of the $\cc=-9$ and $\cc=-3$ (small) $N=4$ superconformal algebra are closely related to other realisations that have appeared in the literature.
\begin{itemize}\itemsep=0pt
	\item For $\smin[1/2]$ ($\cc=-9$), \cite[Proposition~5.1 and Theorem~6.1]{AdaRea14} replaces $\sfvoa$ and $\lvoa$ with the $\cc=-11$ fermionic ghost \svoa\ $\fgvoa$ and the $\cc=2$ bosonic ghost \voa\ $\bgvoa$, respectively.
	We remark that $\sfvoa \subset \fgvoa$ and $\lvoa \supset \bgvoa$ as \svas\ (these are not conformal embeddings).
	\item In \cite[Proposition~1]{AdaVer20}, \smash{$\smin[1/2]$} ($\cc=-9$) is shown to embed into a $\ZZ_2$-orbifold of~\smash{$\sfvoa \otimes \lvoa^{1/2}$}, where~\smash{$\lvoa^{1/2}$} is a simple-current extension of $\lvoa$.
	This is very close to Corollary~\ref{cor:pmhalf}, but is very different in flavour: their realisation uses screening operators and does not identify $\sfvoa$ as a~W-algebra.
	\item On the other hand, \cite[Theorem~4.14]{CreSim15} realises \smash{$\smin[-1/2]$} ($\cc=-3$) as a simple-current extension of the tensor product of the $\cc=-2$ triplet \voa\ $\VOA{W}(2)$ and the $\cc=-1$ affine \voa\ $\ssl[-1/2]$.
	This is also close to Corollary~\ref{cor:pmhalf}: $\VOA{W}(2) \subset \sfvoa$ and~${\ssl[-1/2] \subset \bgvoa \subset \lvoa}$ as \svoas\ (these are conformal embeddings).
\end{itemize}

Recall that $\sfvoa$ has a single irreducible untwisted (hence \ns) module $\sfnsmod \cong \sfvoa$ and a single irreducible $\ZZ_2$-twisted (hence Ramond) module $\sfrmod$.
The \hwvs\ of~$\sfnsmod$ and $\sfrmod$, which we denote by $\sfnsvac$ and $\sfrvac$, have conformal weights $0$ and $-\frac{1}{8}$, respectively.
The Adamovi\'{c} functors thus give two families of almost-irreducible (twisted) \smash{$\sminpmhalf$}-modules, parametrised by $[\lambda] \in \CC/\ZZ$
\begin{equation} \label{eq:defN}
	\nrind{\lambda} = \adfunc{\lambda}\big(\sfnsmod\big) \qquad \text{and} \qquad
	\nnsind{\lambda} = \adfunc{\lambda}\big(\sfrmod\big).
\end{equation}
Note that because \smash{$\lmod{\lambda}$} is $\ZZ_2$-twisted, the \smash{$\nrind{\lambda}$} are too~-- they belong to the Ramond sector of~$\sminpmhalf$.
The \smash{$\nnsind{\lambda}$} are untwisted, hence \ns.

Recall the characters and supercharacters of the irreducible $\sfvoa$-modules
\begin{gather*}
 \fsch{\sfnsmod}{\yy;\qq} = \qq^{1/12} \prod_{n=1}^{\infty} (1-\yy\qq^n)\big(1-\yy^{-1}\qq^n\big)
 = \frac{1}{\yy^{1/2}-\yy^{-1/2}} \frac{\ii \vartheta_{1}(\yy;\qq)}{\eta(\qq)}, \\
 \fch{\sfnsmod}{\yy;\qq} = \qq^{1/12} \prod_{n=1}^{\infty} (1+\yy\qq^n)\big(1+\yy^{-1}\qq^n\big)
 = \frac{1}{\yy^{1/2}+\yy^{-1/2}} \frac{\vartheta_{2}(\yy;\qq)}{\eta(\qq)}, \\
 \fch{\sfrmod}{\yy;\qq} = \qq^{-1/24} \prod_{n=1}^{\infty} \big(1+\yy\qq^{n-1/2}\big)\big(1+\yy^{-1}\qq^{n-1/2}\big)
 = \frac{\vartheta_{3}(\yy;\qq)}{\eta(\qq)}, \\
 \fsch{\sfrmod}{\yy;\qq} = \qq^{-1/24} \prod_{n=1}^{\infty} \big(1-\yy\qq^{n-1/2}\big)\big(1-\yy^{-1}\qq^{n-1/2}\big)
 = \frac{\vartheta_{4}(\yy;\qq)}{\eta(\qq)}.
\end{gather*}
Here, we have included the horizontal grading as in \eqref{eq:uprchar}.
For the lower-bounded $\lvoa$-modules~\smash{$\lmod{\lambda}$}, we define characters as follows
\[
	\fch{\lmod{\lambda}}{\zz;\qq}
	= \traceover{\lmod{\lambda}} \zz^{2b_0} \qq^{t_{0} + (3\kk+2)/12}
	= \frac{\qq^{-(\kk+1)/4 + (3\kk+2)/12}}{\prod_{n=1}^{\infty} (1-\qq^n)^2} \sum_{n \in \ZZ} \zz^{\lambda+n}
	= \zz^{\lambda} \frac{\delta(\zz)}{\eta(\qq)^2}.
\]
Here, $\delta(\zz)$ is the formal series (or distribution) $\sum_{n\in\ZZ} \zz^n$.

Because Corollary~\ref{cor:pmhalf} identifies $J^0$ and $2b$, we may define graded (super)characters of \smash{$\sminpmhalf$}-modules that keep track of the horizontal grading ($\yy$), the $J^0_0$-eigenvalue ($\zz$) and the conformal weight ($\qq$).
For the images \eqref{eq:defN} under the Adamovi\'{c} functors, these (super)characters have the simple forms
\begin{gather}
		\fsch{\nrind{\lambda}}{\yy,\zz;\qq} = \frac{\zz^{\lambda}}{\yy^{1/2}-\yy^{-1/2}} \frac{\ii \vartheta_{1}(\yy;\qq) \delta(\zz)}{\eta(\qq)^3}, \qquad
		\fch{\nnsind{\lambda}}{\yy,\zz;\qq} = \zz^{\lambda} \frac{\vartheta_{3}(\yy;\qq) \delta(\zz)}{\eta(\qq)^3},\nonumber \\
		\fch{\nrind{\lambda}}{\yy,\zz;\qq} = \frac{\zz^{\lambda}}{\yy^{1/2}+\yy^{-1/2}} \frac{\vartheta_{2}(\yy;\qq) \delta(\zz)}{\eta(\qq)^3}, \nonumber \\
		\fsch{\nnsind{\lambda}}{\yy,\zz;\qq} = \zz^{\lambda} \frac{\vartheta_{4}(\yy;\qq) \delta(\zz)}{\eta(\qq)^3}. \label{eq:N=4(s)ch}
\end{gather}
Given these results, we may also use \eqref{eq:N=4sf} to determine the (super)characters of the spectral flows of these modules.

\subsection{Degenerations and spectral flows} \label{sec:degsfk=+-1/2}

Consider now the almost-irreducible \smash{$\sminpmhalf{}$}-modules \smash{$\nrind{\lambda}$}.
These belong to the Ramond sector and so are completely determined by the action of the zero modes, including $G^{\pm}_0$ and $\bG^{\pm}_0$, on their top spaces \cite{DeSFin05}.
A basis for this top space is given by the vectors
\begin{equation} \label{eq:defmuR}
	\rket{\mu} = \sfnsvac \otimes \ee^{-a+\mu c}, \qquad \mu \in [\lambda].
\end{equation}
Using Corollary~\ref{cor:pmhalf}, we compute the action of the generators' zero modes on these top space vectors.
The even modes' actions are given by
\begin{gather}
		J^{+}_{0}\rket{\mu} = \rket{\mu+2}, \qquad
		J^{0}_{0}\rket{\mu} = \mu\rket{\mu}, \qquad
		T_{0}\rket{\mu} = -\tfrac{1}{4}(\kk+1) \rket{\mu}, \nonumber\\
		J^{-}_{0}\rket{\mu} = -\tfrac{1}{4}(\mu - \kk - 1)(\mu + \kk - 1)\rket{\mu-2},
	\qquad \mu \in [\lambda],\label{eq:sl2action}
\end{gather}
while the odd modes all act as $0$.
The top space of \smash{$\nrind{\lambda}$} thus decomposes as a direct sum of two (generically irreducible) modules for the zero-mode algebra (the Ramond-twisted Zhu algebra) of \smash{$\sminpmhalf$}.
We may therefore write
\[ \nrind{\lambda} = \nrrel{\lambda} \oplus \nrrel{\lambda+1},\qquad [\lambda] = \lambda+\ZZ \in \CC/\ZZ,\]
where~${\dbrac{\lambda} = \lambda+2\ZZ \in \CC/2\ZZ}$ and \smash{$\nrrel{\lambda}$} denotes the submodule of \smash{$\nrind{\lambda}$} generated by the \smash{$\rket{\mu}$} with $\mu \in \dbrac{\lambda}$.

Note that the top spaces of the \smash{$\nrrel{\lambda}$}, $\dbrac{\lambda} \in \CC/2\ZZ$, may be viewed as dense modules for the even zero-mode subalgebra of \smash{$\sminpmhalf$} (which is isomorphic to $\gltwo$).
Because the \smash{$\nrrel{\lambda}$} are almost irreducible, see, for example, \cite[Proposition~4.8]{FasMod24}, it follows immediately from \eqref{eq:sl2action} that~\smash{$\nrrel{\lambda}$} is reducible if and only if there exists a \chwv\ $\rket{\mu}$ with~${\mu \in \dbrac{\lambda}}$.
Moreover, such a vector only exists when $\mu = 1\pm\kk = \frac{1}{2}, \frac{3}{2}$.
We thus say that the family $\smash{\bigl\{\nrrel{\lambda} \mid \dbrac{\lambda} \in \CC/2\ZZ}\bigr\}$ degenerates at $\dbrac{\lambda} = \dbrac[\big]{\pm\frac{1}{2}}$.

Let \smash{$\nrihw{\mu}$} denote the irreducible \hw\ \smash{$\sminpmhalf$}-module, in the Ramond sector, whose \hwv\ has $J^0_0$-eigenvalue $\mu$ and $T_0$-eigenvalue $-\frac{1}{4}(\kk+1)$.
Recalling the conjugation automorphism $\nconjsymb$ from Section~\ref{sec:iqhr}, the reducible cases are then characterised by the following nonsplit short exact sequences
\begin{gather}
 \ses{\nconj[\big]{\nrihw{-1/2}}}{\nrrel{1/2}}{\nrihw{-3/2}}, \nonumber\\
 \ses{\nconj[\big]{\nrihw{-3/2}}}{\nrrel{-1/2}}{\nrihw{-1/2}}. \label{ses:N=4R}
\end{gather}
This result deserves a few comments.
\begin{itemize}\itemsep=0pt
	\item The submodules in \eqref{ses:N=4R} are irreducible because any nonzero proper submodule would also be one of \smash{$\nrrel{\pm1/2}$}, hence would intersect the top space.
	\item The quotients are not obviously irreducible, although this would follow easily if one has a~classification of \hwms, such as \cite[Proposition~6.2]{AdaRea14} for $\kk=\frac{1}{2}$.
	Otherwise, the irreducibility may be established using the methods developed in \cite[Proposition~4.13]{KawRel18}, which also demonstrate exactness.
 \item To our knowledge, this is the first time the images of Adamovi{\'c} functors have produced modules that are generically completely reducible rather than generically irreducible. In this case, this may be explained by an observation in \cite{AdaVer20}, whereby \smash{$\sminpmhalf$} in fact embeds into the even subalgebra of $\sfvoa \otimes \lvoa$ with respect to a $\ZZ_2 \times \ZZ_2$ action. The $\ZZ_2$ action on $\sfvoa$ is the usual one for superalgebras. The $\ZZ_2$ action on $\lvoa$ is determined by the action on the lattice generators~$\ee^{mc}$, namely, by whether $m \in \ZZ$ is even or odd. It is clear from inspecting Corollary~\ref{cor:pmhalf} that even generators of \smash{$\sminpmhalf$} are mapped into the $(0,0)$ subspace, and odd generators of \smash{$\sminpmhalf$} are mapped into the $(1,1)$ subspace, with respect to this action.
 \item Finally, the (super)character of \smash{$\nrrel{\lambda}$} may be obtained from that of \smash{$\nrind{\lambda}$}, given in \eqref{eq:N=4(s)ch}, by decomposing into $(\ZZ_2 \times \ZZ_2)$-eigenspaces (where we assign $\rket{\lambda}$ the grade $(0,0)$).
 More precisely, that of \smash{$\nrrel{\lambda}$} is the sum (difference) of the characters of the $(0,0)$- and $(1,1)$-eigenspaces.
\end{itemize}

One can repeat the above analysis in the \ns\ sector.
Now, the top space of \smash{$\nnsind{\lambda}$}, $[\lambda] \in \CC/\ZZ$, has the following basis vectors\[
	\nsket{\mu} = \sfrvac \otimes \ee^{-a+\mu c}, \qquad \mu \in [\lambda].\]
The odd generators have no zero modes, and the even zero modes' actions are given explicitly~by
\begin{gather*}
		J^{+}_{0}\nsket{\mu} = \nsket{\mu+2}, \qquad
		J^{0}_{0}\nsket{\mu} = \mu\nsket{\mu}, \qquad
		T_{0}\nsket{\mu} = -\tfrac{1}{4}\brac[\big]{\kk+\tfrac{3}{2}} \nsket{\mu}, \\
		J^{-}_{0}\nsket{\mu} = -\brac[\big]{\tfrac{1}{4}(\mu - \kk - 1)(\mu + \kk - 1) + \tfrac{1}{8}\bigl(\kk+\tfrac{1}{2}\bigr)(3\kk-1)}\nsket{\mu-2},
	\qquad \mu \in [\lambda].
\end{gather*}
We therefore again get a direct sum decomposition
\[
	\nnsind{\lambda} = \nnsrel{\lambda} \oplus \nnsrel{\lambda+1},\]
where \smash{$\nnsrel{\mu}$} denotes the irreducible \hw\ \smash{$\sminpmhalf$}-module, in the \ns\ sector, whose \hwv\ has $J^0_0$-eigenvalue $\mu$ and $T_0$-eigenvalue $-\frac{1}{4}\big(\kk+\frac{3}{2}\big)$.

\begin{figure}
	\centering
	\begin{tikzpicture}[xscale=1.6]
		\foreach \i in {-3,-1,1,3} \node[wt,label=above:{$\nsket{\i}$}] at (\i,1) {};
		\foreach \i in {-2,0,2} {\node[wt] at (\i,0.1) {}; \node[wt] at (\i,-0.1) {};}
		\node (l) at (-4,1) {$\dots$};
		\node (r) at (4,1) {$\dots$};
		\foreach \i in {-3,3} \node at (\i,0) {$\dots$};
		\foreach \i in {-2,0,2} \node at (\i,-1) {$\vdots$};
		\draw (l) -- (-1,1);
		\draw (1,1) -- (r);
		\node[wt] (-1) at (-1,1) {};
		\node[wt] (+1) at (1,1) {};
		\node[wt] (+0) at (0,0.1) {};
		\node[wt] (-0) at (0,-0.1) {};
		\draw[->,thick,blue] (-1) to[bend left] node[above] {$\scriptstyle J^+_0$} (+1);
		\draw[->,thick,blue] (-1) to[bend left] node[above right, inner sep=0mm] {$\hspace{-3mm}\scriptstyle G^+_{-1/2}$} (+0);
		\draw[->,thick,blue] (+0) to[bend left] node[above left, inner sep=0mm] {$\scriptstyle \bG^+_{1/2}\hspace{-1mm}$} (+1);
		\draw[->,thick,blue] (-1) to[bend right] node[below left, inner sep=0mm] {$\scriptstyle \bG^+_{-1/2}$} (-0);
		\draw[->,thick,blue] (-0) to[bend right] node[below right, inner sep=1mm] {$\scriptstyle G^+_{1/2}$} (+1);
		\node at (5,1) {$\Delta = -\frac{1}{2}$};
		\node at (5,0) {$\Delta = 0$};
	\end{tikzpicture}
\caption{%
	An illustration of the structure of the \smash{$\smin[1/2]$}-module \smash{$\nnsrel{1}$} .
	The dots represent weight vectors, with the $J^0_0$-eigenvalue increasing from left to right and the $T_0$-eigenvalue $\Delta$ from top to bottom.
} \label{fig:nsrel+1/2}
\end{figure}
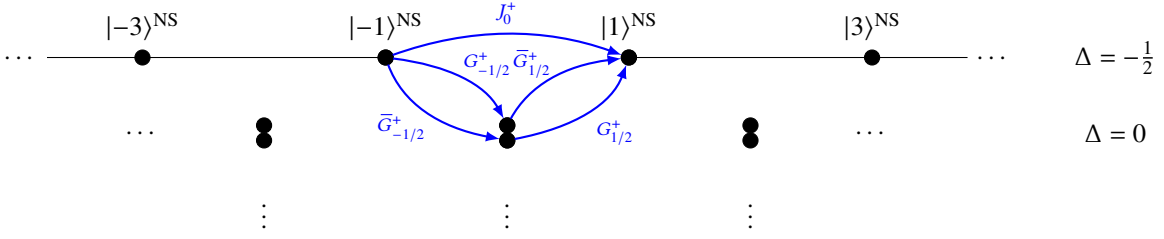

When $\kk=-\frac{1}{2}$, the \ns\ analysis is almost identical to the Ramond case.
The family $\bigl\{\nnsrel{\lambda} \mid \dbrac{\lambda} \in \CC/2\ZZ\bigr\}$ degenerates at $\dbrac{\lambda} = \dbrac{\pm\frac{1}{2}}$ and there we have the following nonsplit short exact sequences:
\begin{gather*} 
 \ses{\nconj[\big]{\nnsihw{-1/2}}}{\nnsrel{1/2}}{\nnsihw{-3/2}}, \\
 \ses{\nconj[\big]{\nnsihw{-3/2}}}{\nnsrel{-1/2}}{\nnsihw{-1/2}}.
\end{gather*}
The analysis for $\kk=\frac{1}{2}$ is significantly more interesting.
In this case, $J^{-}_{0}\nsket{\mu} = 0$ has only one solution: $\mu=1$.
The family $\bigl\{\nnsrel{\lambda} \mid \dbrac{\lambda} \in \CC/2\ZZ\bigr\}$ thus only degenerates at $\dbrac{\lambda} = \dbrac{1}$.
But, \smash{$\nnsrel{1}$} has more than two composition factors!
Indeed, its Loewy diagram is
\[
	\begin{tikzpicture}[scale=0.8,->,baseline=(l.base)]
		\node (t) at (0,2) {$\nnsihw{-1}$};
		\node (l) at (-2,0) {$\nnsihw{0}$};
		\node (r) at (2,0) {$\nnsihw{0}$};
		\node (b) at (0,-2) {$\nconj{\nnsihw{-1}}$};
		\draw (t) -- (l);
		\draw (t) -- (r);
		\draw (l) -- (b);
		\draw (r) -- (b);
		\node[nom] at (0,0) {$\nnsrel{1}$};
	\end{tikzpicture}.
\]
This structure is also indicated schematically in Figure~\ref{fig:nsrel+1/2}.

We make a few further comments on the structure of \smash{$\nnsrel{1}$} .
\begin{itemize}\itemsep=0pt
	\item This module is not logarithmic: $T_0$ acts semisimply (its head and socle are not isomorphic).
	\item It is well known that an almost-irreducible module may have composition factors that are not detected by the top space, see \cite[Section~4.5]{KawRel18}.
	For simple affine \vas\ (and W-algebras), this appears to be a feature of the representation theory at nonadmissible levels.
	Given that the level of the affine \vsa\ generated by $J^0$ and $J^{\pm}$ is admissible for~${\kk=-\frac{1}{2}}$, but nonadmissible for~${\kk=\frac{1}{2}}$, this qualitative structural difference should perhaps not be surprising.
	\item It is nevertheless quite satisfying to use \eqref{eq:simplequotemb} to explicitly verify that \smash{$G^-_{-1/2}$} and \smash{$\bG^-_{-1/2}$} annihilate \smash{$\nsket{1}$}, while \smash{$G^-_{1/2}$} and \smash{$\bG^-_{1/2}$} annihilate both \smash{$G^+_{-1/2} \nsket{-1}$} and \smash{$\bG^+_{-1/2} \nsket{-1}$}.
	\item The fact that \smash{$\nnsrel{1}$} does not have more than four composition factors uses the classification \cite[Proposition~6.2]{AdaRea14} of irreducible \hw\ \smash{$\smin[1/2]$}-modules \big(there are only two: \smash{$\nnsihw{0}$} and \smash{$\nnsihw{-1}$}\big).
	The inverse-reduction technology developed in \cite{AdaWei23} can also be used to prove this classification result, though we shall not do so here.
	\item As noted above, \smash{$\smin[1/2]$} is a \vsa\ of the tensor product of a $\cc=-11$ fermionic ghost system and a $\cc=2$ bosonic ghost system.
	In our language, \cite[Proposition~6.3]{AdaRea14} identifies this tensor product as a nonsplit extension of $\nnsihw{-1}$ by $\nnsihw{0}$, hence as a length-$2$ quotient of~\smash{$\nnsrel{1}$}.\looseness=1
	\item The quotient \smash{$\nnsrel{1} \big/ \nconj[\big]{\nnsihw{-1}}$} is a \hw\ \smash{$\smin[1/2]$}-module with two linearly independent singular vectors that share the same $J^0_0$- and $T^0$-eigenvalues.
	Both generate a copy of the vacuum module $\nnsihw{0}$.
\end{itemize}

We conclude by reporting the action of spectral flow on the irreducible \hw\ \smash{$\sminpmhalf$}-modules.
For $\kk=-\frac{1}{2}$, the spectral flow orbits feature ordinary modules that do not appear as composition factors of the almost-irreducible modules constructed by inverse reduction.
These follow from explicit computation using \eqref{eq:N=4sf}
\begin{gather*}
\cdots \overset{\nsfsymb{1/2}}{\lra} \nconj[\big]{\nrihw{-1/2}} \overset{\nsfsymb{1/2}}{\lra} \nnsihw{0} \overset{\nsfsymb{1/2}}{\lra} \nrihw{-1/2} \overset{\nsfsymb{1/2}}{\lra} \cdots, \\
		\cdots \overset{\nsfsymb{1/2}}{\lra} \nconj[\big]{\nnsihw{-1/2}} \overset{\nsfsymb{1/2}}{\lra} \nrihw{0} \overset{\nsfsymb{1/2}}{\lra} \nnsihw{-1/2} \overset{\nsfsymb{1/2}}{\lra} \cdots, \\
		\cdots \overset{\nsfsymb{1/2}}{\lra} \nconj[\big]{\nrihw{-3/2}} \overset{\nsfsymb{1/2}}{\lra} \nnsihw{1} \overset{\nsfsymb{1/2}}{\lra} \nrihw{-3/2} \overset{\nsfsymb{1/2}}{\lra} \cdots, \\
		\cdots \overset{\nsfsymb{1/2}}{\lra} \nconj[\big]{\nnsihw{-3/2}} \overset{\nsfsymb{1/2}}{\lra} \nrihw{1} \overset{\nsfsymb{1/2}}{\lra} \nnsihw{-3/2} \overset{\nsfsymb{1/2}}{\lra} \cdots.
\end{gather*}
Here, the dots indicate irreducible modules that are not lower bounded.
We also note the similarity to the spectral flow orbits of $\ssl[-1/2]$, see \cite[Figure~3]{RidSL208}.
Contrarily, for $\kk=\frac{1}{2}$, the orbits contain no new lower-bounded modules
\begin{gather*}
		\cdots \overset{\nsfsymb{1/2}}{\lra} \nconj[\big]{\nrihw{-3/2}} \overset{\nsfsymb{1/2}}{\lra} \nnsihw{0} \overset{\nsfsymb{1/2}}{\lra} \nrihw{-3/2} \overset{\nsfsymb{1/2}}{\lra} \cdots, \\
		\cdots \overset{\nsfsymb{1/2}}{\lra} \nconj[\big]{\nrihw{-1/2}} \overset{\nsfsymb{1/2}}{\lra} \nnsihw{-1} \overset{\nsfsymb{1/2}}{\lra} \cdots, \qquad
		\cdots \overset{\nsfsymb{1/2}}{\lra} \nconj[\big]{\nnsihw{-1}} \overset{\nsfsymb{1/2}}{\lra} \nrihw{-1/2} \overset{\nsfsymb{1/2}}{\lra} \cdots.
	\end{gather*}

\subsection[Logarithmic sminpmhalf-modules]{Logarithmic \smash{$\boldsymbol{\sminpmhalf}$}-modules} \label{sec:k=+-1/2log}

Recall that the Ramond sector of the symplectic fermions \svoa\ $\sfvoa$ is semisimple, but the \ns\ sector is not.
In particular, it contains a logarithmic module~$\sfnsmod_{(4)}$.
This is induced from the projective module of \smash{$\psl$} and is thus an indecomposable sum of four copies of the vacuum module.
We present its Loewy diagram and a basis for its top space\looseness=-1
\begin{equation} \label{eq:sfproj}
	\begin{tikzpicture}[scale=0.8,->,baseline=(l.base)]
		\node (t) at (0,2) {$\sfnsmod$};
		\node (l) at (-2,0) {$\sfnsmod$};
		\node (r) at (2,0) {$\sfnsmod$};
		\node (b) at (0,-2) {$\sfnsmod$};
		\draw (t) -- (l);
		\draw (t) -- (r);
		\draw (l) -- (b);
		\draw (r) -- (b);
		\node[nom] at (0,0) {$\sfnsmod_{(4)}$};
		\begin{scope}[shift={(8,0)}]
			\node (t) at (0,2) {$\ket{\NS;t}$};
			\node (l) at (-2,0) {$\ket{\NS;m}$};
			\node (r) at (2,0) {$\ket{\NS;\bm}$};
			\node (b) at (0,-2) {$\ket{\NS;b}$};
			\draw (t) -- node[above left] {$\chi_0$} (l);
			\draw (t) -- node[above right] {$\bchi_0$} (r);
			\draw (l) -- node[below left] {$\bchi_0$} (b);
			\draw (r) -- node[below right] {$-\chi_0$} (b);
		\end{scope}
	\end{tikzpicture}
	.
\end{equation}
The logarithmic nature of \smash{$\sfnsmod_{(4)}$} is evident from the action of the Virasoro zero mode: it maps $\ket{\NS;t}$ to $\frac{1}{2} \ket{\NS;b}$.

We can apply Adamovi\'{c} functors to such nonsemisimple $\sfvoa$-modules.
To warm up, we will first consider a structurally simpler example: the length-$2$ submodule of \smash{$\sfnsmod_{(4)}$} generated by~$\ket{\NS;m}$.
This $\sfvoa$-module, which we shall denote by $\sfnsmod_{(2)}$, has top space $\spn \{\ket{\NS;m}, \ket{\NS;b}\}$.
It is moreover a self-extension of the vacuum module that coincides with a Verma module of~$\apsl$.
We shall denote its images under the Adamovi\'{c} functors by $\nrverma{\lambda}$, $[\lambda] \in \CC/\ZZ$.
These are lower-bounded Ramond-twisted \smash{$\sminpmhalf$}-modules.

The top space of $\nrverma{\lambda}$ has basis $\bigl\{\rket{\mu;m}, \rket{\mu;b} \mid \mu \in [\lambda]\bigr\}$, where we follow \eqref{eq:defmuR} but add the label ``$m$'' or ``$b$'' appropriately.
Using Corollary~\ref{cor:pmhalf}, it is easy to determine the action of the zero-mode algebra of \smash{$\sminpmhalf$} on this basis.
Indeed, the even zero modes act as in \eqref{eq:sl2action} (preserving labels), and almost all of the odd zero modes act as zero.
The exceptions are those involving~$\bchi_0$
\[ 
	\bG^+_0 \rket{\mu;m} = \rket{\mu+1;b} \qquad \text{and} \qquad
	\bG^-_0 \rket{\mu;m} = -\tfrac{1}{2}\bigl(\mu-\tfrac{1}{2}\bigr) \rket{\mu-1;b}.
\]
Because of this, we still have a direct sum decomposition
\smash{$\nrverma{\lambda} = \nrver{\lambda} \oplus \nrver{\lambda+1}$}, \smash{$ [\lambda] \in \CC/\ZZ$},
where~\smash{$\nrver{\lambda}$} is spanned by the $\rket{\mu;m}$ and $\rket{\mu+1;b}$ with $\mu \in \dbrac{\lambda}$.
However, the \smash{$\nrver{\lambda}$} are nonsemisimple and (generically) have two composition factors
\begin{equation} \label{eq:sesV}
	\ses{\nrrel{\lambda+1}}{\nrver{\lambda}}{\nrrel{\lambda}}, \qquad \dbrac{\lambda} \ne \dbrac[\big]{\pm\tfrac{1}{2}}.
\end{equation}

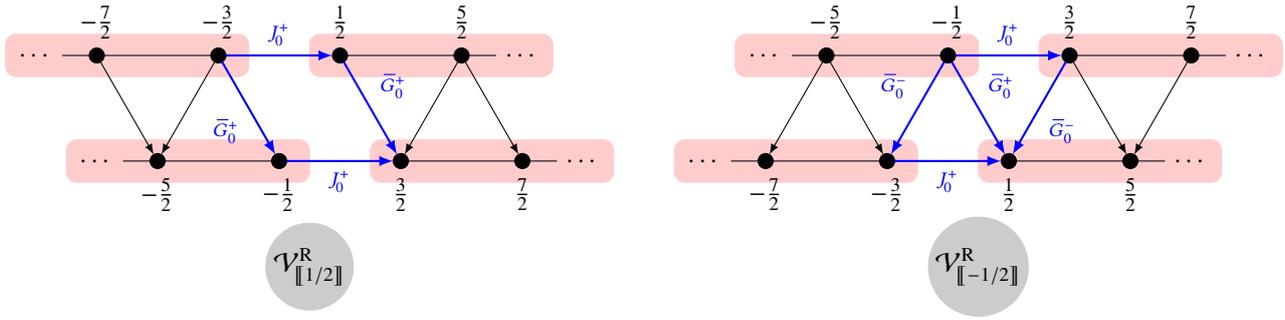
\begin{figure}
	\centering
	\begin{tikzpicture}[xscale=0.75,yscale=1.4]
		\filldraw[red!20!white,rounded corners] (-5,1.2) rectangle (-1,0.8);
		\filldraw[red!20!white,rounded corners] (-4,0.2) rectangle (0,-0.2);
		\filldraw[red!20!white,rounded corners] (4,1.2) rectangle (0,0.8);
		\filldraw[red!20!white,rounded corners] (5,0.2) rectangle (1,-0.2);
		\node[wt,label=above:{$-\frac{7}{2}$}] (-72) at (-7/2,1) {};
		\node[wt,label=above:{$-\frac{3}{2}$}] (-32) at (-3/2,1) {};
		\node[wt,label=above:{$\frac{1}{2}$}] (+12) at (1/2,1) {};
		\node[wt,label=above:{$\frac{5}{2}$}] (+52) at (5/2,1) {};
		\node[wt,label=below:{$\frac{7}{2}$}] (+72) at (7/2,0) {};
		\node[wt,label=below:{$\frac{3}{2}$}] (+32) at (3/2,0) {};
		\node[wt,label=below:{$-\frac{1}{2}$}] (-12) at (-1/2,0) {};
		\node[wt,label=below:{$-\frac{5}{2}$}] (-52) at (-5/2,0) {};
		\node (tl) at (-9/2,1) {$\dots$};
		\node (tr) at (7/2,1) {$\dots$};
		\node (br) at (9/2,0) {$\dots$};
		\node (bl) at (-7/2,0) {$\dots$};
		\draw (tl) -- (-32)
		 (tr) -- (+12)
		 (bl) -- (-12)
		 (br) -- (+32);
		\draw[->] (-72) -- (-52);
		\draw[->] (-32) -- (-52);
		\draw[->]	(+52) -- (+32);
		\draw[->] (+52) -- (+72);
		\draw[->,thick,blue] (-32) -- node[above] {$\scriptstyle J^+_0$} (+12);
		\draw[->,thick,blue] (-12) -- node[below] {$\scriptstyle J^+_0$} (+32);
		\draw[->,thick,blue] (-32) -- node[below left] {$\scriptstyle \bG^+_0$} (-12);
		\draw[->,thick,blue] (+12) -- node[above right] {$\scriptstyle \bG^+_0$} (+32);
		\node[nom] at (0,-1) {$\nrver{1/2}$};
		\begin{scope}[shift={(11,0)}]
			\filldraw[red!20!white,rounded corners] (-5,0.2) rectangle (-1,-0.2);
			\filldraw[red!20!white,rounded corners] (-4,1.2) rectangle (0,0.8);
			\filldraw[red!20!white,rounded corners] (4,0.2) rectangle (0,-0.2);
			\filldraw[red!20!white,rounded corners] (5,1.2) rectangle (1,0.8);
			\node[wt,label=below:{$-\frac{7}{2}$}] (-72) at (-7/2,0) {};
			\node[wt,label=below:{$-\frac{3}{2}$}] (-32) at (-3/2,0) {};
			\node[wt,label=below:{$\frac{1}{2}$}] (+12) at (1/2,0) {};
			\node[wt,label=below:{$\frac{5}{2}$}] (+52) at (5/2,0) {};
			\node[wt,label=above:{$\frac{7}{2}$}] (+72) at (7/2,1) {};
			\node[wt,label=above:{$\frac{3}{2}$}] (+32) at (3/2,1) {};
			\node[wt,label=above:{$-\frac{1}{2}$}] (-12) at (-1/2,1) {};
			\node[wt,label=above:{$-\frac{5}{2}$}] (-52) at (-5/2,1) {};
			\node (tl) at (-9/2,0) {$\dots$};
			\node (tr) at (7/2,0) {$\dots$};
			\node (br) at (9/2,1) {$\dots$};
			\node (bl) at (-7/2,1) {$\dots$};
			\draw (tl) -- (-32)
			 (tr) -- (+12)
			 (bl) -- (-12)
			 (br) -- (+32);
			\draw[<-] (-72) -- (-52);
			\draw[<-] (-32) -- (-52);
			\draw[<-]	(+52) -- (+32);
			\draw[<-] (+52) -- (+72);
			\draw[->,thick,blue] (-32) -- node[below] {$\scriptstyle J^+_0$} (+12);
			\draw[->,thick,blue] (-12) -- node[above] {$\scriptstyle J^+_0$} (+32);
			\draw[<-,thick,blue] (-32) -- node[above left] {$\scriptstyle \bG^-_0$} (-12);
			\draw[<-,thick,blue] (+12) -- node[below right] {$\scriptstyle \bG^-_0$} (+32);
			\draw[->,thick,blue] (-12) -- node[above right] {$\scriptstyle \bG^+_0$} (+12);
			\node[nom] at (0,-1) {$\nrver{-1/2}$};
		\end{scope}
	\end{tikzpicture}
\caption{%
	An illustration of the structure of the top spaces of the \smash{$\sminpmhalf$}-modules $\nrver{1/2}$ (left) and~$\nrver{-1/2}$ (right).
	The dots represent weight vectors, with the $J^0_0$-eigenvalue increasing from left to right.
	The top row corresponds to vectors with label ``$m$'' and the bottom row has label ``$b$''.
	The shading suggests the composition factors.
} \label{fig:Vpm1/2}
\end{figure}

When $\dbrac{\lambda} = \dbrac[\big]{\pm\frac{1}{2}}$, there is a further degeneration corresponding to the fact that \smash{$\nrrel{\pm1/2}$} is reducible.
In this case, $J^-_0$ annihilates $\rket{\mu;\bullet}$, for $\mu=\frac{1}{2}, \frac{3}{2}$ and $\bullet=m, b$, and $\bG^-_0$ annihilates~$\rket{\frac{1}{2};m}$.
The resulting Loewy diagrams are then
\begin{equation} \label{eq:Vpm1/2}
	\begin{tikzpicture}[scale=0.8,->,baseline=(l.base)]
		\node (t) at (0,2) {$\nrihw{-3/2}$};
		\node (l) at (-2,0) {$\nrihw{-1/2}$};
		\node (r) at (2,0) {$\nconj[\big]{\nrihw{-1/2}}$};
		\node (b) at (0,-2) {$\nconj[\big]{\nrihw{-3/2}}$};
		\draw (t) -- (l);
		\draw (t) -- (r);
		\draw (l) -- (b);
		\draw (r) -- (b);
		\node[nom] at (0,0) {$\nrver{1/2}$};
		\begin{scope}[shift={(8,0)}]
			\node (t) at (0,2) {$\nrihw{-1/2}$};
			\node (l) at (-2,0) {$\nrihw{-3/2}$};
			\node (r) at (2,0) {$\nconj[\big]{\nrihw{-3/2}}$};
			\node (b) at (0,-2) {$\nconj[\big]{\nrihw{-1/2}}$};
			\draw (t) -- (l);
			\draw (t) -- (r);
			\draw (l) -- (b);
			\draw (r) -- (b);
			\node[nom] at (0,0) {$\nrver{-1/2}$};
		\end{scope}
	\end{tikzpicture}
	\ ,
\end{equation}
see also Figure~\ref{fig:Vpm1/2}.
Of course, neither module is logarithmic.

With this analysis complete, we turn to the image of the logarithmic module \smash{$\sfnsmod_{(4)}$} under the Adamovi{\'c} functor, which we denote by \smash{$\nrlog{\lambda}$}, $[\lambda] \in \CC/\ZZ$.
The top space has basis $\bigl\{ \rket{\mu;t}, \rket{\mu;m}, \rket{\mu;\bm}, \rket{\mu;b} \mid \mu \in [\lambda] \bigr\}$, where again the second label follows~\eqref{eq:sfproj}.
On this top space, the even zero modes $J^{+}_{0}$ and $J^{0}_{0}$ act as they did on $\nrverma{\lambda}$ (preserving labels).
The remaining even zero modes act, for $\bullet = t,m,\bm,b$, as follows
\begin{gather}
		J^{-}_{0} \rket{\mu;\bullet} = -\tfrac{1}{4}(\mu-\kk-1)(\mu+\kk-1) \rket{\mu-2;\bullet} + \tfrac{1}{2}\bigl(\kk+\tfrac{1}{2}\bigr)(3\kk-1) \rket{\mu-2;b} \delta_{\bullet,t}, \nonumber\\
		T_0 \rket{\mu;\bullet} = -\tfrac{1}{4}(\kk+1)\rket{\mu;\bullet} + \tfrac{1}{2} \rket{\mu;b} \delta_{\bullet,t}.\label{eq:sl2action''}
\end{gather}
The action of the odd zero modes on the top space basis vectors is again zero, except for
\begin{gather*} 
 G^{+}_{0} \rket{\mu;t} = \rket{\mu+1;m}, \qquad G^{+}_{0} \rket{\mu;\bm} = -\rket{\mu+1;b}, \\
 \bG^{+}_{0} \rket{\mu;t} = \rket{\mu+1;\bm}, \qquad \bG^{+}_{0} \rket{\mu;m} = +\rket{\mu+1;b}, \\
 G^{-}_{0} \rket{\mu;t} = -\tfrac{1}{2}\bigl(\mu-\tfrac{1}{2}\bigr) \rket{\mu-1;m},\qquad
 G^{-}_{0} \rket{\mu;\bm} = +\tfrac{1}{2}\bigl(\mu-\tfrac{1}{2}\bigr) \rket{\mu-1;b}, \\
 \bG^{-}_{0} \rket{\mu;t} = -\tfrac{1}{2}\bigl(\mu-\tfrac{1}{2}\bigr) \rket{\mu-1;\bm},\qquad
 \bG^{-}_{0} \rket{\mu;m} = -\tfrac{1}{2}\bigl(\mu-\tfrac{1}{2}\bigr) \rket{\mu-1;b}.
\end{gather*}

It is clear that we once again have a direct sum decomposition
\[
\nrlog{\lambda} = \nrsublog{\lambda} \oplus \nrsublog{\lambda+1} , \qquad [\lambda] \in \CC/\ZZ,
\]
where $\nrsublog{\lambda}$ is spanned by the $\rket{\mu;t}$, $\rket{\mu+1;m}$, $\rket{\mu+1;\bm}$ and $\rket{\mu;b}$, with $\mu \in \dbrac{\lambda}$.
Each of these submodules is logarithmic and their Loewy diagrams have the generic form
\[
	\begin{tikzpicture}[scale=0.8,->,baseline=(l.base)]
		\node (t) at (0,2) {$\nrrel{\lambda}$};
		\node (l) at (-2,0) {$\nrrel{\lambda+1}$};
		\node (r) at (2,0) {$\nrrel{\lambda+1}$};
		\node (b) at (0,-2) {$\nrrel{\lambda}$};
		\draw (t) -- (l);
		\draw (t) -- (r);
		\draw (l) -- (b);
		\draw (r) -- (b);
		\node[nom] at (0,0) {$\nrsublog{\lambda}$};
	\end{tikzpicture}
	, \qquad \dbrac{\lambda} \ne \dbrac[\big]{\pm\tfrac{1}{2}}.
\]
It is clear from \eqref{eq:sl2action''} that these modules are logarithmic.

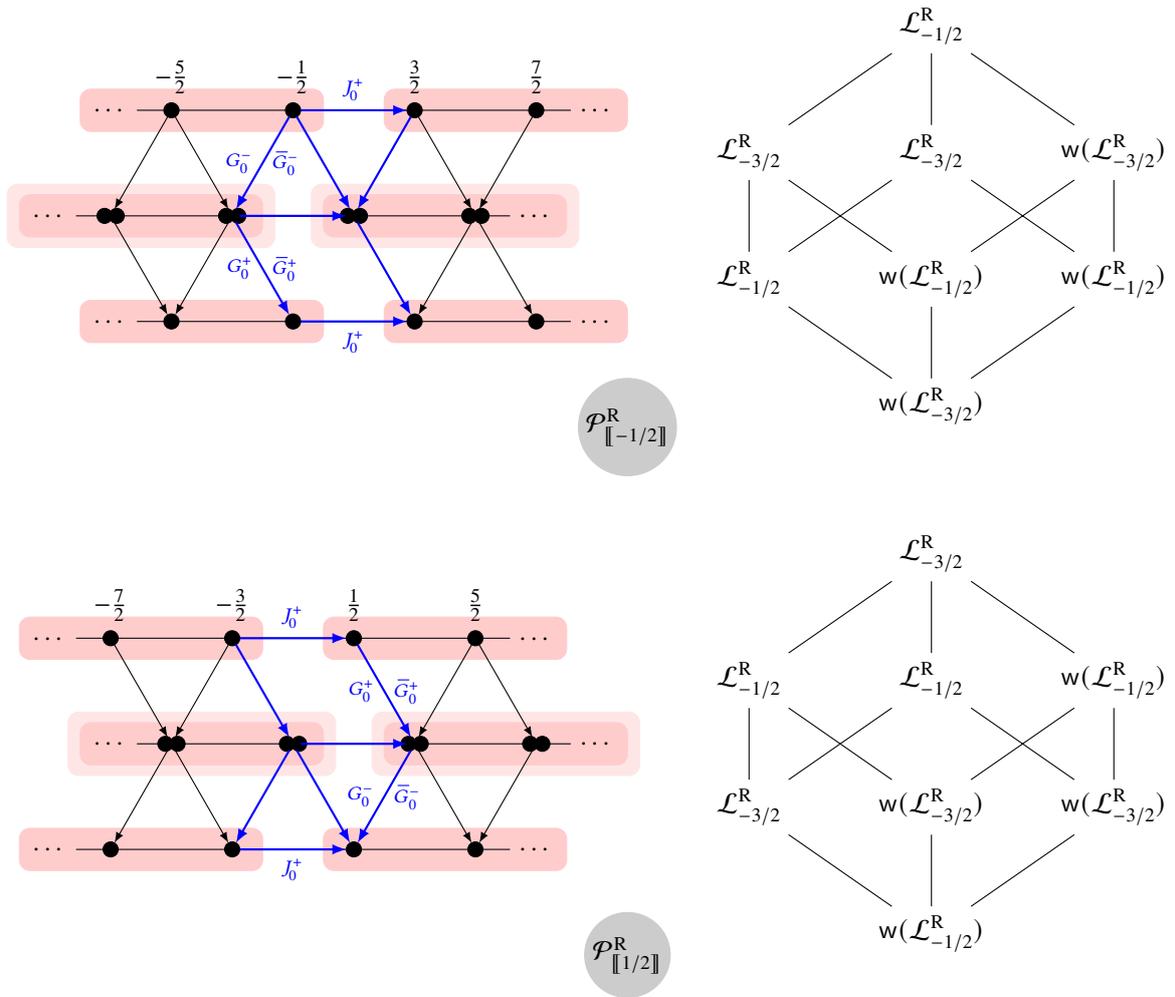
\begin{figure}[t]
	\centering
	\begin{tikzpicture}[xscale=0.8,yscale=1.4]
		\filldraw[red!20!white,rounded corners] (1,0.2) rectangle (5,-0.2);
		\filldraw[red!10!white,rounded corners] (-0.2,1.3) rectangle (4.2,0.7);
		\filldraw[red!20!white,rounded corners] (0,1.2) rectangle (4,0.8);
		\filldraw[red!20!white,rounded corners] (1,2.2) rectangle (5,1.8);
		\filldraw[red!20!white,rounded corners] (0,0.2) rectangle (-4,-0.2);
		\filldraw[red!10!white,rounded corners] (-0.8,1.3) rectangle (-5.2,0.7);
		\filldraw[red!20!white,rounded corners] (-1,1.2) rectangle (-5,0.8);
		\filldraw[red!20!white,rounded corners] (0,2.2) rectangle (-4,1.8);
		\node[wt,label=above:{$\frac{7}{2}$}] (+72t) at (7/2,2) {}; 
		\node[wt,label=above:{$\frac{3}{2}$}] (+32t) at (3/2,2) {};
		\node[wt,label=above:{$-\frac{1}{2}$}] (-12t) at (-1/2,2) {};
		\node[wt,label=above:{$-\frac{5}{2}$}] (-52t) at (-5/2,2) {};
		\node[wt,opacity=0] (-72m) at (-7/2,1) {}; 
		\node[wt,opacity=0] (-32m) at (-3/2,1) {};
		\node[wt,opacity=0] (+12m) at (1/2,1) {};
		\node[wt,opacity=0] (+52m) at (5/2,1) {};
		\node[wt] at (-7/2-0.1,1) {}; 
		\node[wt] at (-3/2-0.1,1) {};
		\node[wt] at (1/2-0.1,1) {};
		\node[wt] at (5/2-0.1,1) {};
		\node[wt] at (-7/2+0.1,1) {}; 
		\node[wt] at (-3/2+0.1,1) {};
		\node[wt] at (1/2+0.1,1) {};
		\node[wt] at (5/2+0.1,1) {};
		\node[wt] (+72b) at (7/2,0) {}; 
		\node[wt] (+32b) at (3/2,0) {};
		\node[wt] (-12b) at (-1/2,0) {};
		\node[wt] (-52b) at (-5/2,0) {};
		\node (tr) at (9/2,2) {$\dots$};
		\node (tl) at (-7/2,2) {$\dots$};
		\node (ml) at (-9/2,1) {$\dots$};
		\node (mr) at (7/2,1) {$\dots$};
		\node (br) at (9/2,0) {$\dots$};
		\node (bl) at (-7/2,0) {$\dots$};
		\draw (tl) -- (tr) 
		 (ml) -- (mr)
		 (bl) -- (br);
		\draw[->] (-52t) -- (-32m); 
		\draw[->] (+32t) -- (+52m);
		\draw[->] (-72m) -- (-52b);
		\draw[->] (+52m) -- (+72b);
		\draw[->] (-52t) -- (-72m); 
		\draw[->] (+72t) -- (+52m);
		\draw[->] (-32m) -- (-52b);
		\draw[->] (+52m) -- (+32b);
		\draw[->,thick,blue] (-12t) -- node[above] {$\scriptstyle J^+_0$} (+32t);
		\draw[->,thick,blue] (-32m) -- (+12m);
		\draw[->,thick,blue] (-12b) -- node[below] {$\scriptstyle J^+_0$} (+32b);
		\draw[->,thick,blue] (-12t) -- (+12m);
		\draw[->,thick,blue] (+12m) -- (+32b);
		\draw[->,thick,blue] (-32m) -- node[left] {$\scriptstyle G^+_0$} node[right] {$\scriptstyle \bG^+_0$} (-12b);
		\draw[->,thick,blue] (-12t) -- node[left] {$\scriptstyle G^-_0$} node[right] {$\scriptstyle \bG^-_0$} (-32m);
		\draw[->,thick,blue] (+32t) -- (+12m);
		\node[nom] at (5,-1) {$\nrsublog{-1/2}$};
		\begin{scope}[shift={(10,1)},yscale=0.6]
			\node (1) at (0,3) {$\nrihw{-1/2}$};
			\node (21) at (-3,1) {$\nrihw{-3/2}$};
			\node (22) at (0,1) {$\nrihw{-3/2}$};
			\node (23) at (3,1) {$\nconj[\big]{\nrihw{-3/2}}$};
			\node (31) at (-3,-1) {$\nrihw{-1/2}$};
			\node (32) at (0,-1) {$\nconj[\big]{\nrihw{-1/2}}$};
			\node (33) at (3,-1) {$\nconj[\big]{\nrihw{-1/2}}$};
			\node (4) at (0,-3) {$\nconj[\big]{\nrihw{-3/2}}$};
			\draw (1) -- (21)
			 (1) -- (22)
			 (1) -- (23)
			 (21) -- (31)
			 (21) -- (32)
			 (22) -- (31)
			 (22) -- (33)
			 (23) -- (32)
			 (23) -- (33)
			 (31) -- (4)
			 (32) -- (4)
			 (33) -- (4);
		\end{scope}
		\begin{scope}[shift={(0,-5)}]
			\filldraw[red!20!white,rounded corners] (0,0.2) rectangle (4,-0.2);
			\filldraw[red!10!white,rounded corners] (0.8,1.3) rectangle (5.2,0.7);
			\filldraw[red!20!white,rounded corners] (1,1.2) rectangle (5,0.8);
			\filldraw[red!20!white,rounded corners] (0,2.2) rectangle (4,1.8);
			\filldraw[red!20!white,rounded corners] (-1,0.2) rectangle (-5,-0.2);
			\filldraw[red!10!white,rounded corners] (0.2,1.3) rectangle (-4.2,0.7);
			\filldraw[red!20!white,rounded corners] (0,1.2) rectangle (-4,0.8);
			\filldraw[red!20!white,rounded corners] (-1,2.2) rectangle (-5,1.8);
			\node[wt,label=above:{$\frac{5}{2}$}] (+52t) at (5/2,2) {}; 
			\node[wt,label=above:{$\frac{1}{2}$}] (+12t) at (1/2,2) {};
			\node[wt,label=above:{$-\frac{3}{2}$}] (-32t) at (-3/2,2) {};
			\node[wt,label=above:{$-\frac{7}{2}$}] (-72t) at (-7/2,2) {};
			\node[wt,opacity=0] (-52m) at (-5/2,1) {}; 
			\node[wt,opacity=0] (-12m) at (-1/2,1) {};
			\node[wt,opacity=0] (+32m) at (3/2,1) {};
			\node[wt,opacity=0] (+72m) at (7/2,1) {};
			\node[wt] at (-5/2-0.1,1) {}; 
			\node[wt] at (-1/2-0.1,1) {};
			\node[wt] at (3/2-0.1,1) {};
			\node[wt] at (7/2-0.1,1) {};
			\node[wt] at (-5/2+0.1,1) {}; 
			\node[wt] at (-1/2+0.1,1) {};
			\node[wt] at (3/2+0.1,1) {};
			\node[wt] at (7/2+0.1,1) {};
			\node[wt] (+52b) at (5/2,0) {}; 
			\node[wt] (+12b) at (1/2,0) {};
			\node[wt] (-32b) at (-3/2,0) {};
			\node[wt] (-72b) at (-7/2,0) {};
			\node (tr) at (7/2,2) {$\dots$};
			\node (tl) at (-9/2,2) {$\dots$};
			\node (ml) at (-7/2,1) {$\dots$};
			\node (mr) at (9/2,1) {$\dots$};
			\node (br) at (7/2,0) {$\dots$};
			\node (bl) at (-9/2,0) {$\dots$};
			\draw (tl) -- (tr) 
			 (ml) -- (mr)
			 (bl) -- (br);
			\draw[->] (-72t) -- (-52m); 
			\draw[->] (+52t) -- (+72m);
			\draw[->] (-52m) -- (-32b);
			\draw[->] (+32m) -- (+52b);
			\draw[->] (-32t) -- (-52m); 
			\draw[->] (+52t) -- (+32m);
			\draw[->] (-52m) -- (-72b);
			\draw[->] (+72m) -- (+52b);
			\draw[->,thick,blue] (-32t) -- node[above] {$\scriptstyle J^+_0$} (+12t);
			\draw[->,thick,blue] (-12m) -- (+32m);
			\draw[->,thick,blue] (-32b) -- node[below] {$\scriptstyle J^+_0$} (+12b);
			\draw[->,thick,blue] (-32t) -- (-12m);
			\draw[->,thick,blue] (-12m) -- (+12b);
			\draw[->,thick,blue] (+12t) -- node[left] {$\scriptstyle G^+_0$} node[right] {$\scriptstyle \bG^+_0$} (+32m);
			\draw[->,thick,blue] (-12m) -- (-32b);
			\draw[->,thick,blue] (+32m) -- node[left] {$\scriptstyle G^-_0$} node[right] {$\scriptstyle \bG^-_0$} (+12b);
			\node[nom] at (5,-1) {$\nrsublog{1/2}$};
			\begin{scope}[shift={(10,1)},yscale=0.6]
				\node (1) at (0,3) {$\nrihw{-3/2}$};
				\node (21) at (-3,1) {$\nrihw{-1/2}$};
				\node (22) at (0,1) {$\nrihw{-1/2}$};
				\node (23) at (3,1) {$\nconj[\big]{\nrihw{-1/2}}$};
				\node (31) at (-3,-1) {$\nrihw{-3/2}$};
				\node (32) at (0,-1) {$\nconj[\big]{\nrihw{-3/2}}$};
				\node (33) at (3,-1) {$\nconj[\big]{\nrihw{-3/2}}$};
				\node (4) at (0,-3) {$\nconj[\big]{\nrihw{-1/2}}$};
				\draw (1) -- (21)
				 (1) -- (22)
				 (1) -- (23)
				 (21) -- (31)
				 (21) -- (32)
				 (22) -- (31)
				 (22) -- (33)
				 (23) -- (32)
				 (23) -- (33)
				 (31) -- (4)
				 (32) -- (4)
				 (33) -- (4);
			\end{scope}
		\end{scope}
	\end{tikzpicture}
\caption{%
	Loewy diagrams and top-space structures of the modules $\nrsublog{\pm1/2}$ when $\kk = -\frac{1}{2}$.
	The top row of the structures corresponds to vectors with label ``$t$'', the middle row to those with labels ``$m$'' and~``$\protect\bm$'', and the bottom row to those with labels ``$b$''.
	The shading again suggests the composition factors.
} \label{fig:P1/2}
\end{figure}

Consider the degeneration that occurs when $\dbrac{\lambda} = \dbrac[\big]{\pm\frac{1}{2}}$.
We will suppose first that $\kk=-\frac{1}{2}$.
Then, the qualitative difference from the generic action given above is that $J^-_0$ acts as zero, when $\mu=\frac{1}{2}$ or $\frac{3}{2}$, while $G^-_0$ and $\bG^-_0$ act as zero, when $\mu=\frac{1}{2}$.
The resulting structures are depicted, along with the corresponding Loewy diagrams, in Figure~\ref{fig:P1/2}.

For $\kk=\frac{1}{2}$, the only difference is that $J^-_0$ no longer acts as zero when $\mu=\frac{1}{2}$ or $\frac{3}{2}$.
Instead, we have
\[
	J^-_0 \rket{\mu;\bullet} = \tfrac{1}{4} \rket{\mu-2;b} \delta_{\bullet,t}.
\]
This modifies the structural depictions and Loewy diagrams of Figure~\ref{fig:P1/2} by adding one arrow to each.
Specifically, we should add an arrow in each depiction from the top dot labelled by~$\frac{1}{2}$~\big($\frac{3}{2}$\big) to the bottom dot labelled by~$-\frac{3}{2}$~\big($-\frac{1}{2}$\big) and an arrow in each Loewy diagram from the second-row factor $\nconj[\big]{\nrihw{-1/2}}$ \big($\nconj[\big]{\nrihw{-3/2}}$\big) to the third-row factor $\nrihw{-3/2}$ \big($\nrihw{-1/2}$\big).

We finish with three remarks.
\begin{itemize}\itemsep=0pt
	\item The structures of the \smash{$\nrsublog{\lambda}$}, for both $\kk=\frac{1}{2}$ and $-\frac{1}{2}$, may be partially characterised in terms of short exact sequences involving the \smash{$\nrver{\lambda}$} of \eqref{eq:sesV} and \eqref{eq:Vpm1/2}
	\[
		\ses{\nrver{\lambda+1}}{\nrsublog{\lambda}}{\nrver{\lambda}}.
	\]
	This unsurprising characterisation is nevertheless not as informative as the Loewy diagrams of Figure~\ref{fig:P1/2}.
	\item Given experience with other ``logarithmic'' \svoas, see \cite{CreLog13} for example, one might be surprised by the existence of the family of logarithmic \smash{$\sminpmhalf$}-modules \smash{$\bigl\{\nrsublog{\lambda} \mid \dbrac{\lambda} \in \CC/2\ZZ\bigr\}$}.
	Indeed, the standard module formalism of \cite{CreLog13,RidVer14} suggests that the ``typical'' members of the family of ``standard modules'' should be both irreducible and projective, hence logarithmic modules should be (in this sense) rare.
	In the case at hand, the ubiquity of logarithmic modules is a simple consequence of the fact that \smash{$\spr[\pm1/2]$} is not rational and, in fact, admits a logarithmic module (see \cite[Section~4.2]{AdaWei23} for another example of this phenomenon).
	Examples like these indicate that the simple ``typical/atypical'' distinction of the standard module formalism needs generalising.
	\item A natural question is whether the \smash{$\nrsublog{\lambda}$} are projective in a suitable category of \smash{$\sminpmhalf$}-modules.
	Experience with projectives for $\ssl$ \cite{AraWei23,CreCos18} and (conjecturally) $\saff{-3/2}{\slthree}$ \cite{CreKaz21} suggests that the answer is no because the \smash{$\nrsublog{\lambda}$} are lower bounded.
	It would be extremely interesting to investigate the construction of logarithmic \smash{$\sminpmhalf$}-modules that do not have a lower-bounded spectral flow image.
	We expect that these modules will feature Jordan blocks of rank greater than $2$ in the action of $T_0$.
\end{itemize}

\subsection*{Acknowledgements}

We thank Dra\v{z}en Adamovi\'{c}, Thomas Creutzig, Justine Fasquel, Andy Linshaw and Sven M{\"o}ller for helpful conversations on topics related to this research.
In particular, we thank Dra\v{z}en Adamovi\'{c}, for suggesting a proof of the simplicity of the principal W-algebra (see Theorem~\ref{thm:prsimple}), and one of the anonymous reviewers, for sketching another interesting alternative proof of the same fact.
ZF's research is supported by the Australian Research Council Discovery Project DP210101502.
CR's research is supported by the Deutsche Forschungsgemeinschaft through the Collaborative Research Centre 1624: Higher Structures, Moduli Spaces and Integrability, project number 50663264.
DR's research is supported by the Australian Research Council Discovery Project DP210101502 and an Australian Research Council Future Fellowship FT200100431.

\pdfbookmark[1]{References}{ref}
\LastPageEnding

\end{document}